\documentclass[11pt]{amsart}
\usepackage{amsxtra}
\usepackage{amssymb}
\usepackage[mathscr]{eucal}
\theoremstyle{plain}
\usepackage{color}

\newcommand{\nn}{\hfill\nonumber}
\newtheorem{theorem}{Theorem}[subsection]
\newtheorem{lemma}[theorem]{Lemma}
\newtheorem{definition-theorem}[theorem]{Definition-Theorem}
\newtheorem{proposition}[theorem]{Proposition}
\newtheorem{corollary}[theorem]{Corollary}
\newtheorem{definition}[theorem]{Definition}
\newtheorem{example}[theorem]{Example}
\newtheorem{remark}[theorem]{Remark}
\newtheorem{conjecture}[theorem]{Conjecture}
\newtheorem{notation}[theorem]{Notation}

\newtheorem*{maintheorem*}{Main Theorem}
\newcommand \bth[1] { \begin{theorem}\label{t#1} }
\newcommand \ble[1] { \begin{lemma}\label{l#1} }

\newcommand \bpr[1] { \begin{proposition}\label{p#1} }
\newcommand \bco[1] { \begin{corollary}\label{c#1} }
\newcommand \bde[1] { \begin{definition}\label{d#1}\rm }
\newcommand \bex[1] { \begin{example}\label{e#1}\rm }
\newcommand \bre[1] { \begin{remark}\label{r#1}\rm }
\newcommand \bcj[1] { \begin{conjecture}\label{j#1}\rm }

\newcommand \bnota[1] { \begin{notation}\label{n#1}\rm }
\renewcommand {\eth} { \end{theorem} }
\newcommand {\ele} { \end{lemma} }

\newcommand {\epr} { \end{proposition} }
\newcommand {\eco} { \end{corollary} }
\newcommand {\ede} { \end{definition} }
\newcommand {\eex} { \end{example} }
\newcommand {\ere} { \end{remark} }
\newcommand {\ecj} { \end{conjecture} }

\newcommand {\enota} { \end{notation} }
\newcommand \thref[1]{Theorem \ref{t#1}}
\newcommand \leref[1]{Lemma \ref{l#1}}
\newcommand \prref[1]{Proposition \ref{p#1}}
\newcommand \coref[1]{Corollary \ref{c#1}}

\newcommand \exref[1]{Example \ref{e#1}}

\makeatletter
\newsavebox{\@brx}
\newcommand{\llangle}[1][]{\savebox{\@brx}{\(\m@th{#1\langle}\)}%
  \mathopen{\copy\@brx\kern-0.5\wd\@brx\usebox{\@brx}}}
\newcommand{\rrangle}[1][]{\savebox{\@brx}{\(\m@th{#1\rangle}\)}%
  \mathclose{\copy\@brx\kern-0.5\wd\@brx\usebox{\@brx}}}
\makeatother

\usepackage{tikz}
\usetikzlibrary{matrix}
\usepackage{graphicx,ctable,booktabs}
\usetikzlibrary{arrows.meta}
\usepackage{tikz-cd}

\DeclareMathOperator{\Ind}{Ind}

\DeclareMathOperator{\Ext}{Ext} 
\DeclareMathOperator{\hHom}{\widehat{Hom}} 
\DeclareMathOperator{\hEnd}{\widehat{End}}

\DeclareMathOperator{\Spec}{Spec} 
\DeclareMathOperator{\Spech}{Spec^h} 
\DeclareMathOperator{\Ann}{Ann}

 \DeclareMathOperator{\Proj}{Proj}

\DeclareMathOperator{\Aut}{Aut}

\DeclareMathOperator{\End}{End} 
 \DeclareMathOperator{\Hom}{Hom}

\DeclareMathOperator{\ThickId}{ThickId}
\DeclareMathOperator{\Mod}{{\sf Mod}}
\DeclareMathOperator{\modd}{{\sf mod}}

\DeclareMathOperator{\Spc}{Spc}

\DeclareMathOperator{\ev}{ev}
\DeclareMathOperator{\coev}{coev}

\DeclareMathOperator{\op}{op}

\DeclareMathOperator{\cone}{cone}

\DeclareMathOperator{\Set}{{\sf Set}}

\DeclareMathOperator{\opH}{\operatorname{H}}

\newcommand{\mf}{\mathfrak}
\newcommand{\mc}{\mathcal}
\newcommand{\mb}{\mathbb}
\newcommand{\id}{\operatorname{id}}

\newcommand{\kk}{\Bbbk}
\newcommand{\bT}{\mathbf T}

\newcommand{\bK}{\mathbf K}

\newcommand{\bP}{\mathbf P}
\newcommand{\bQ}{\mathbf Q}
\newcommand{\bZ}{\mathbf Z}
\newcommand{\bI}{\mathbf I}
\newcommand{\bJ}{\mathbf J}

\newcommand{\Loc}{\operatorname{Loc}}
\newcommand{\XX}{\mathcal X}

\newcommand{\unit}{\ensuremath{\mathbf 1}}

\newcounter{listequation}

\numberwithin{equation}{section}

\begin{document}
\title[On the Spectrum and Support Theory of 
a Finite Tensor Category]
{On the Spectrum and Support Theory of \\
a Finite Tensor Category}
\author[Daniel K. Nakano]{Daniel K. Nakano}
\address{Department of Mathematics \\
University of Georgia \\
Athens, GA 30602\\
U.S.A.}
\thanks{Research of D.K.N. was supported in part by NSF grants DMS-1701768 and DMS-2101941.}
\email{nakano@math.uga.edu}

\author[Kent B. Vashaw]{Kent B. Vashaw}
\address{
Department of Mathematics\\
Massachusetts Institute of Technology\\
Cambridge, MA 02139\\
U.S.A.}
\thanks{Research of K.B.V. was supported in part by a Board of Regents LSU fellowship, an Arthur K. Barton Superior Graduate Student Scholarship in Mathematics from LSU, NSF grant DMS-1901830 and NSF Postdoctoral Fellowship DMS-2103272.}
\email{kentv@mit.edu}
\author[Milen T. Yakimov]{Milen T. Yakimov}
\address{
Department of Mathematics \\
Northeastern University, Boston \\
MA 02115 \\
U.S.A.}
\thanks{Research of M.T.Y. was supported in part by NSF grants DMS-2131243 and  DMS-2200762.}
\email{m.yakimov@northeastern.edu}
\subjclass[2010]{Primary 18M05, 18G65; Secondary 16T05}

\keywords{finite tensor categories, tensor triangular geometry, Hopf algebras}

\begin{abstract}
Finite tensor categories (FTCs) $\bT$ are important generalizations of the categories of finite dimensional modules of 
finite dimensional Hopf algebras, which play a key role in many areas of mathematics and mathematical physics. 
There are two fundamentally different support theories for them: a cohomological one and a universal one 
based on the noncommutative Balmer spectra of their stable (triangulated) categories $\underline \bT$. 

In this paper we introduce the key notion of the categorical center $C^\bullet_{\underline \bT}$ of the cohomology ring $R^\bullet_{\underline \bT}$ of an FTC, $\bT$. 
This enables us to put forward a complete and detailed program to investigate the relationship between the two support theories, based 
on $C^\bullet_{\underline \bT}$ of the cohomology ring $R^\bullet_{\underline \bT}$ of an FTC, $\bT$. 

Our main result is the construction of a continuous map from the noncommutative Balmer spectrum of an arbitrary FTC, $\bT$, 
to the $\Proj$ of the categorical center $C^\bullet_{\underline \bT}$ and a theorem that this map is surjective under 
a weaker finite generation assumption for $\bT$ than the one conjectured by Etingof-Ostrik. 
We conjecture that, for all FTCs, (i) the map is a homeomorphism and (ii) the two-sided thick ideals of $\underline \bT$
are classified by the specialization closed subsets of $\Proj C^\bullet_{\underline \bT}$. We verify parts of the conjecture 
under stronger assumptions on the category $\bT$. 

Many examples are presented that demonstrate how in important cases $C^\bullet_{\underline \bT}$ arises as a fixed point 
subring of $R^\bullet_{\underline \bT}$ and how the two-sided thick ideals of $\underline \bT$ are determined in a uniform fashion
(while previous methods dealt on a case-by-case basis with case specific methods).
The majority of our results are proved in the greater generality of monoidal triangulated categories
and versions of them for Tate cohomology are also presented.
\end{abstract}
\maketitle

\section{Introduction} \label{Intro}

\subsection{\ } 
\label{1.1}
Fifty years ago, Quillen pioneered the idea that even though the cohomology ring for a finite group cannot be completely calculated, its spectrum can be described via elementary abelian 
subgroups \cite{Quillen1}. This promoted the spectrum as an important object in the cohomology and representation theory of finite groups.  The theory of support varieties was developed
in the 1980's to bridge homological properties of modules to the ambient geometry of the cohomological spectrum. In a major breakthrough, Friedlander and Suslin in 1995 proved that the cohomology ring for a finite group scheme (equivalently, a finite dimensional cocommutative Hopf algebra) is finitely generated \cite{FS1}, and the theory of support varieties was developed in this setting through the innovation of $\pi$-points by 
Friedlander and Pevtsova \cite{FP1}. 

The subject evolved further by the mid 2000's when Balmer introduced tensor triangular geometry (TTG) to provide a unifying method for addressing problems in representation theory,
homotopy theory and algebraic geometry \cite{Balmer1, Balmer2}. Given a symmetric monoidal category, $\bf K$, Balmer's idea was to utilize the monoidal (tensor) structure 
to view the category like a commutative ring, and construct a topological space $\Spc \bK$ 
via prime ideals in $\bf K$. Many earlier ideas in the cohomology theory fit nicely into the framework of TTG. For instance, in the representations of finite group 
schemes and quantum groups at roots of unity, the topological spaces $\Spc \bK$ manifest themselves naturally as the projectivization of the spectrum of the cohomology ring.

The authors \cite{NVY1} developed a general noncommutative version of Balmer's theory that deals with an arbitrary monoidal triangulated 
category $\bK$ (M$\Delta$C for short).  There are many important families of M$\Delta$Cs coming from diverse areas of mathematics and mathematical physics. 
One such family is comprised of the stable module categories of finite dimensional Hopf algebras, and more generally, the stable categories of finite tensor categories. 
These M$\Delta$Cs are in general not symmetric (because Hopf algebras are not necessarily cocommutative). 
The key feature of our new approach is to define the noncommutative Balmer spectrum $\Spc \bK$ and support data for $\bK$
in terms of tensoring thick ideals of $\bK$, and not to use object-wise tensoring. This enables us to distinguish between prime ideals 
and completely prime ideals in $\Spc \bK$, and provide new insights about the tensor product property on supports. 

\subsection{\ } 
\label{1.2}
In recent years, there has been a substantial effort to better understand the cohomology and support theory for finite dimensional Hopf algebras $H$, and more generally, 
for finite tensor categories $\bT$. By definition, a finite tensor category is an abelian $k$-linear category with $k$-bilinear monoidal product, such that morphism spaces are finite-dimensional, all objects have finite length, the unit object is simple, there are enough projectives, there are finitely many isomorphism classes of simple objects, and every object has a left and right dual. Equivalently, $\bT$ is a $k$-linear monoidal category which is equivalent, as an abelian category, to the category of modules for a finite-dimensional $k$-algebra; see \cite[Defs. 1.8.5 and 4.1.1 and Ch. 6]{EGNO}
for background. Finite tensor categories are important generalizations of the categories of finite dimensional modules $\modd(H)$
and play an important role in mathematical physics and quantum computing. 

There are two major approaches to support theory for finite tensor categories:
\begin{enumerate}
\item[(i)] cohomological support and
\item[(ii)] support via the (noncommutative) Balmer spectrum of the stable category $\underline{\bT}$.
\end{enumerate}
For the second one, we note that all finite tensor categories are Frobenius \cite{EO}. The corresponding stable categories $\underline{\bT}$
are M$\Delta$Cs, which have the property that they are the compact parts of M$\Delta$Cs admitting 
arbitrary set indexed coproducts.

In this paper we provide a program with the goal of establishing a precise relationship between the two support theories 
for all finite tensor categories. Many parts of this program apply to the compact subcategories of 
all M$\Delta$Cs admitting arbitrary set indexed coproducts.

\subsection{\ } 
\label{1.3}
For a topological space $X$, denote by $\mc{X}_{cl}(X)$, the collection of all of its closed subsets. 
Define the {\em{cohomology ring}} of an M$\Delta$C, $\bK$ to be the ring 
\[
\quad R^\bullet_{\bK} :=\bigoplus_{i \geq 0} \Hom_{\bK}(\unit, \Sigma^i \unit),
\]
where $\unit$ denotes the unit object of $\bK$. Define the {\em{Tate cohomology ring}} of $\bK$ to be the ring
\[
\widehat R^\bullet_{\bK}:= \bigoplus_{i \in \mathbb{Z}} \Hom_{\bK}(\unit, \Sigma^i \unit).
\]
It is well known that these rings are (graded) commutative (see e.g. \cite{SA}).
For $M,N \in \bK$ denote the sets
\[
\Hom^\bullet_{\bK} ( M, N) = \bigoplus_{i \geq 0} \Hom_{\bK}(M, \Sigma^i N), \quad
\hHom^\bullet_{\bK} ( M, N) = \bigoplus_{i \in \mathbb{Z}} \Hom_{\bK}(M, \Sigma^i N).
\]
They have canonical structures of $R^\bullet_{\bK}$ and $\widehat R^\bullet_{\bK}$-modules, respectively, given by
\[
g \cdot h = \Sigma^i(h) (g \otimes \id_M), 
\quad \mbox{for} \; \; g \in \Hom_{\bK}(\unit, \Sigma^i \unit), h \in \Hom_{\bK}(M, \Sigma^j N),
\]
where we use the natural isomorphism $(\Sigma^i \unit) \otimes M \cong \Sigma^i M$. 
Set $\End^\bullet_{\bK} (M) = \Hom^\bullet_{\bK} ( M, M)$, 
$\hEnd^\bullet_{\bK} (M) = \hHom^\bullet_{\bK} ( M, M)$.

The {\em{cohomological support map}} $W : \bK \to \mc{X}_{cl} (\Proj  R_{\bK}^\bullet)$ and 
{\em{the Tate cohomological support map}} $\widehat W : \bK \to \mc{X}_{cl} (\Proj  \widehat R_{\bK}^\bullet)$ are given by
\begin{align*}
W(M) &= \{ \mf{p} \in \Proj R^\bullet_{\bK} : {\mf{p}} \supseteq \Ann_{R^\bullet_{\bK}} \big( \End^\bullet_{\bK} (M) \big)  \},
\\
\widehat W(M) &= \{ \mf{p} \in \Spech \widehat R^\bullet_{\bK} : {\mf{p}} \supseteq \Ann_{\widehat R^\bullet_{\bK}} \big( \hEnd^\bullet_{\bK} (M) \big) \}.
\end{align*}
Here, following \cite{Balmer1}, for a $\mb Z$-graded commutative ring $R$ we denote by 
\[
\Spech R
\]
the space of homogeneous prime ideals of $R$, equipped with the Zariski topology.
For a commutative ring $R$ and a finitely generated $R$-module $Y$, 
\[
\{ {\mf{p}} \in \Spec R : Y_{\mf{p}} \neq 0 \} =  \{ {\mf{p}} \in \Spec R : {\mf{p}} \supseteq \Ann_R(Y) \}, 
\]
see e.g. \cite[Ch.\ 3, Exercise 19(v)]{AM}. Therefore, the cohomological support map and the Tate cohomological support map
are given by the more classical formulas 
\begin{align*}
W(M) &= \{ \mf{p} \in \Proj R^\bullet_{\bK} : \End^\bullet_{\bK} (M)_{\mf{p}} \not = 0 \}, \\
\widehat W(M) &= \{ \mf{p} \in \Spech \widehat R^\bullet_{\bK} : \hEnd^\bullet_{\bK} (M)_{\mf{p}} \not = 0 \},
\end{align*}
respectively  for those objects $M \in \bK$ for which $\End^\bullet_{\bK} ( M)$ is a finitely generated $R^\bullet_{\bK}$-module and 
$\hEnd^\bullet_{\bK} ( M)$ is a finitely generated $\widehat R^\bullet_{\bK}$-module.
These supports maps were widely used for many classes of M$\Delta$Cs in various degrees of generality 
(stable module categories of finite groups, finite group schemes, and finite dimensional Hopf algebras, derived categories of algebraic varieties, etc.),  
and were defined in this generality in \cite{BKSS}.

The noncommutative Balmer spectrum $\Spc \bK$ of $\bK$ is the topological 
space of all thick prime ideals of $\bK$, see Section \ref{noncommBalmerspec} for definitions. 
The basic notions were introduced by Buan, Krause and Solberg \cite{BKS}. The concepts were reintroduced and extensively studied 
by the authors in \cite{NVY1}.
The corresponding support map $V : \bK \to \mc{X}_{cl} (\Spc \bK)$ is defined by
\[
V(M) := \{ \bP \in \Spc \bK \mid M \notin \bP \}.
\]
The following is one of the major problems in the area:
\medskip

\noindent
{\bf{Problem A.}} {\em{What is the relationship between the cohomological support map 
$W : \bK \to \mc{X}_{cl} (\Proj R_{\bK}^\bullet)$ and the Tate cohomological support map $\widehat W : \bK \to \mc{X}_{cl} (\Spech  \widehat R_{\bK}^\bullet)$, on the one hand, 
and the Balmer spectrum support map $V : \bK \to \mc{X}_{cl} (\Spc \bK)$, on the other hand?
}}
\noindent

\subsection{\ } 
\label{1.4}
We now present a new framework with the aim of settling this problem for the stable categories of all finite tensor categories. 

Let ${\mathcal K}$ be a collection of objects that generates the M$\Delta$C, $\bK$, as a triangulated category. Define the {\em{categorical center}} $\widehat C_{\bK}^\bullet$ of the Tate cohomology ring $\widehat R_{\bK}^\bullet$ as the graded subring spanned by all $g \in \Hom_{\bK}(\unit, \Sigma^i \unit)$, such that for every object $M\in {\mathcal K}$, the following diagram commutes, where the isomorphisms are structure isomorphisms for an M$\Delta$C:
\begin{equation}
\label{diag-M-g}
\begin{tikzcd}
\unit \otimes M \arrow[d, "g \otimes \id_M"] \arrow[r, "\cong"] & M          & M \otimes \unit \arrow[l, "\cong"] \arrow[d, "\id_M \otimes g"] \\
\Sigma^i \unit \otimes M \arrow[r, "\cong"]                     & \Sigma^i M & M \otimes \Sigma^i \unit \arrow[l, "\cong"]                    
\end{tikzcd}
\end{equation}
Define the categorical center of the 
cohomology ring $R_{\bK}^\bullet$ to be
\[
C_{\bK}^\bullet := R_{\bK}^\bullet \cap \widehat C_{\bK}^\bullet.
\]

When $\mc{K}$ is taken to be the collection of all objects of $\bK$, then $\widehat C_{\bK}^\bullet$ is isomorphic to the Tate cohomology ring of the Drinfeld center of $\bK$ (cf. \cite[Section 7.3]{EGNO}), 
which is a monoidal additive category with an autoequivalence induced from $\Sigma$. 
For smaller collections $\mc{K}$ (which matter in our situation),  
our categorical centers $\widehat C_{\bK}^\bullet$ are generally bigger algebras than those. 
(Note that the Tate cohomology ring $\widehat  R_{\bK}^\bullet$ is itself graded commutative for all M$\Delta$Cs, $\bK$, 
so $\widehat C_{\bK}^\bullet$ is not a center in the classical sense). 
If $\bK$ is a braided M$\Delta$C, then $\widehat C_\bK^\bullet = \widehat R_\bK^\bullet$. 
The categorical center $\widehat C_{\bK}^\bullet$ of the Tate cohomology ring $\widehat R_{\bK}^\bullet$ 
is a subalgebra of the 
graded center of $\bK$ studied in \cite{BF,BIK1} and other papers,
but the latter is typically much bigger and contains $\widehat R_{\bK}^\bullet$.
In Section \ref{sec:five-algebras} we present a detailed 
comparison of the categorical center with other algebras of homological origin. 
We refer the reader to the diagram \eqref{eq:nine-algs} that summarizes this relationship.
\medskip
\\
\noindent
{\bf{Remark.}} In both notations $\widehat C_{\bK}^\bullet$ and  $C_{\bK}^\bullet$ we suppress the dependence on a generating set of objects $\mc{K}$
to avoid overcrowded notation.
\medskip

In Section \ref{sec:fixed-pt} we describe how in important situations, 
the categorical centers $C_\bK^\bullet$ and $\widehat C_\bK^\bullet$ arise as 
fixed point subring of $R_\bK^\bullet$ and $\widehat R_\bK^\bullet$ under group actions. 
In Section \ref{sec:examples} we present an explicit calculation of the 
categorical centers $C_\bK^\bullet$ and $\widehat C_\bK^\bullet$ 
for the stable module categories of many classes of 
(non-cocommutative) Hopf algebras that have been actively studied before.
\medskip
\\
\noindent
{\bf{Remark.}} In the case when $\bK$ is the stable category $\underline{\bT}$ of
a finite tensor category $\bT$, the collection ${\mathcal K}$ will be chosen to be the (finite) set 
of simple objects of $\bT$.
\medskip

With the notion of categorical center, we present the following key definition that will be used throughout this paper. 
\medskip
\\
\noindent
{\bf{Definition.}}
The {\em{central cohomological support}} and {\em{Tate central cohomological support}} of an M$\Delta$C, $\bK$, are the maps 
$W_C : \bK \to \Proj C_{\bK}^\bullet$ and $\widehat W_C : \bK \to \Proj \widehat C_{\bK}^\bullet$ given by
\begin{align*}
W_C(M) &= \{ \mf{p} \in \Proj C^\bullet_{\bK} : {\mf{p}} \supseteq \Ann_{C^\bullet_{\bK}} \big( \End^\bullet_{\bK} (M) \big)  \},
\\
\widehat W_C(M) &= \{ \mf{p} \in \Spech \widehat C^\bullet_{\bK} : {\mf{p}} \supseteq \Ann_{\widehat C^\bullet_{\bK}} \big( \hEnd^\bullet_{\bK} (M) \big) \}.
\end{align*}

\subsection{\ } 
\label{1.5}
Etingof and Ostrik \cite{EO} conjectured that all finite tensor categories $\bT$ satisfy the strong finite generation condition:
\vskip .25cm 
{\em (fg) $R_{\underline{\bT}}^\bullet$ is noetherian and for all $M \in \underline{\bT}$, $\End^\bullet_{\underline{\bT}} (M)$ is a finitely generated $R_{\underline{\bT}}^\bullet$-module.}
\vskip .25cm 
\noindent 
We say that an M$\Delta$C, $\bK$, satisfies the {\em{weak finite generation condition}} if 
\smallskip

(wfg) for all $M \in \bK$, $\End^\bullet_{\bK} (M)$ is a finitely generated $C_{\bK}^\bullet$-module
\smallskip
\\ \noindent
and that $\bK$ satisfies the {\em{weak Tate finite generation condition}} if 
\smallskip

(wTfg) for all $M \in \bK$, $\hEnd^\bullet_{\bK} (M)$ is a finitely generated $\widehat C_{\bK}^\bullet$-module. 
\smallskip
\\ \noindent
Since for all $M, N \in \bK$, $\Hom^\bullet_{\bK} (M,N)$ is a direct summand of $\End^\bullet_{\bK} (M \oplus N)$
and $\hHom^\bullet_{\bK} (M,N)$ is a direct summand of $\hEnd^\bullet_{\bK} (M \oplus N)$, the two conditions are equivalent to 
\smallskip

(wfg) for all $M, N \in \bK$, $\Hom^\bullet_{\bK} (M,N)$ is a finitely generated $C_{\bK}^\bullet$-module and

(wTfg) for all $M, N \in \bK$, $\hHom^\bullet_{\bK} (M, N)$ is a finitely generated $\widehat C_{\bK}^\bullet$-module. 
\smallskip
\\ \noindent
Here, the terminology ``weak" is justified by the fact that if the Drinfeld center of a finite tensor category $\bT$ satisfies the (fg) condition, then $\underline{\bT}$ satisfies the (wfg) condition; hence, the validity of the Etingof--Ostrik conjecture implies that all stable categories of finite tensor categories satisfy (wfg). To see this, we first note that Negron and Plavnik \cite{NP1} proved that
the validity of the Etingof--Ostrik conjecture implies that the cohomology ring of $\bT$ is finitely generated over the cohomology ring of the Drinfeld center of $\bT$. The first algebra is isomorphic to $R^\bullet_{\underline{\bT}}$ and the action of the cohomology ring of the Drinfeld center of $\bT$ factors through $C^\bullet_{\underline{\bT}}$, 
see Section \ref{sec:five-algebras} and 
\cite[Proposition 3.3]{NP1}.

Our main result proves the existence of continuous maps $\rho$ and $\widehat \rho$ that can be used to solve Problem A for all FTCs.
It is striking that such continuous maps can be constructed without any assumptions on the M$\Delta$C, $\bK$, see part (a) of 
the theorem. In parts (b)-(d) we establish properties of these maps under weak finite generation assumptions.
\medskip

\noindent
{\bf{Theorem B.}} {\em{Let $\bK$ be an M$\Delta$C.}}
\begin{enumerate} 
\item[(a)] {\em{There is a continuous map $\widehat \rho: \Spc \bK \to \Spech \widehat  C^\bullet_{\bK}$ defined by 
\begin{align*}
\widehat \rho: \Spc \bK &\to \Spech \widehat C^\bullet_{\bK} \\
\bP &\mapsto \langle g \in \widehat C^\bullet_{\bK}: \cone (g) \not \in \bP \rangle.
\end{align*}
Its composition with the contraction map $\Spech \widehat C^\bullet_{\bK} \to \Spech C^\bullet_{\bK}$, 
gives rise to the continuous map 
\[
\rho: \Spc \bK \to \Spech C^\bullet_{\bK}, \quad
\bP \mapsto \langle g \in C^\bullet_{\bK}: \cone (g) \not \in \bP \rangle.
\]
In addition, we have the compatibility properties
\[
V(A) \subseteq \widehat  \rho^{-1} (\widehat W_C(A)) 
\quad \mbox{and} \quad
V(A) \subseteq \rho^{-1} (W_C(A)) 
\]
for all $A \in \bK$.}}
\item[(b)] {\em{If $\bK$ satisfies the weak finite generation condition, then the image of $\rho$ contains $\Proj C^\bullet_{\bK}$}}.
\item[(c)] {\em{If $\bK$ is the stable category of a finite tensor category and satisfies the weak finite generation condition 
then $\rho$ takes values in $\Proj C^\bullet_{\bK}$ and is surjective.}}
\item[(d)] {\em{If $\bK$ satisfies the weak Tate finite generation condition, then $\widehat \rho$ is surjective}}.
\end{enumerate}
\medskip

The first part of the theorem is motivated by a result of Balmer \cite[Theorem 5.3]{Balmer2} which established its validity in the symmetric monoidal case.
In connection to part (c) of the theorem, 
in the symmetric case, Balmer also proved that $\widehat \rho$ is surjective under the 
assumption that $\bK$ is connective (i.e. that $R^\bullet_{\bK} = \widehat R^\bullet_{\bK}$), \cite[Theorem 7.3]{Balmer2}, or under the assumption that $\widehat R^\bullet_{\bK}$
is coherent, \cite[Theorems 7.13]{Balmer2}. Both assumptions are strong even in the case of the stable module category of a finite group. 

On the other hand, part (b) of the theorem has important general applications.
As discussed above, the Etingof--Ostrik conjecture implies (wfg) for the stable category of each finite tensor category $\underline{\bT}$; that is, Theorem B(b) and the  Etingof--Ostrik conjecture
imply the surjectivity of $\rho$ for all such categories $\underline{\bT}$. 

Part (c) of the theorem provides a continuous surjective map 
\[
 \Spc (\underline{\modd}(A)) \to \Proj \opH^\bullet(A,k)
\]
for all quasitriangular finite dimensional Hopf algebras satisfying the weak finite generation condition. This is a very general fact, 
a special case of which handles quantum groups of all types.  
Even in the case of the small quantum group of type $A$, the surjectivity is non-trivial, and was independently proved by Negron and Pevtsova in
\cite{NP2}.

\subsection{\ } 
\label{1.6}
We conjecture that the surjective continuous map $\rho$ in Theorem B is a homeomorphism for every finite tensor category (Conjecture E below).
In the rest of the paper we obtain auxiliary results towards the validity of the conjecture under stronger assumptions
(Theorems C and D below). 
Recall from \cite{BKS,NVY1} that a map $\sigma : \bK \to \mc{X}_{cl} (X)$ for a topological space $X$ 
is a {\em{weak support datum map}} if the following conditions are satisified:
\begin{enumerate}
\item $\sigma(0)= \varnothing$ and $\sigma(\unit)= X$;
\item $\sigma(A\oplus B)=\sigma(A)\cup \sigma(B)$, $\forall A, B \in \bK$;
\item $\sigma( \Sigma A) = \sigma(A)$ for all $A \in \bK$;
\item If $A \to B \to C \to \Sigma A$ is a distinguished triangle, then $\sigma(A) \subseteq \sigma(B) \cup \sigma(C)$;
\item $\Phi_\sigma(\bI \otimes \bJ) = \Phi_\sigma(\bI) \cap \Phi_\sigma(\bJ)$ for all thick ideals $\bI$ and $\bJ$ of $\bK$, 
where \\ $\Phi_\sigma(\bI) := \cup_{A \in \bI} \sigma(A)$.
\end{enumerate}
\medskip

We note that for every M$\Delta$C, $\bK$, the Balmer spectrum support map $V : \bK \to \mc{X}_{cl} (\Spc \bK)$ is a weak support datum map. Furthermore, 
the cohomological support map $W : \bK \to \mc{X}_{cl} (\Proj R_{\bK}^\bullet)$, the Tate cohomological support map $\widehat W : \bK \to \mc{X}_{cl} (\Spech \widehat R_{\bK}^\bullet)$
and their central counterparts $W_C : \bK \to \mc{X}_{cl} (\Proj C_{\bK}^\bullet)$ and $\widehat W_C : \bK \to \mc{X}_{cl} (\Spech \widehat C_{\bK}^\bullet)$
always satisfy conditions (1--4), but for arbitrary M$\Delta$Cs it is unknown whether they satisfy condition (5). 
If the central cohomological support map $W_C : \bK \to \mc{X}_{cl} (\Proj C_{\bK}^\bullet)$, 
respectively the Tate central cohomological support map $\widehat W_C : \bK \to \mc{X}_{cl} (\Spech \widehat C_{\bK}^\bullet)$, are 
weak support data maps, then by \cite[Theorem 4.5.1]{NVY1}, we get universal continuous maps
\[
\eta : \Proj C^\bullet_\bK \to \Spc \bK
\quad \mbox{respectively} \quad
\widehat \eta : \Spech \widehat C^\bullet_\bK \to \Spc \bK.
\]

Our first auxiliary theorem constructs right inverses of the surjective continuous maps $\rho$ and $\widehat \rho$ from the main Theorem B 
under certain stronger assumptions on the underlying category.
\medskip

\noindent
{\bf{Theorem C.}} 
\begin{enumerate}
\item[(a)] {\em{If $\bK$ is an M$\Delta$C for which the central cohomological support is a weak support datum and $\bK$ satisfies the (wfg) condition, 
then $\rho( \eta( \mf{p})) = \mf{p}$ for every homogeneous prime ideal $\mf{p}$ of $C^\bullet_\bK$.}}
\item[(b)] {\em{If $\bT$ is a finite tensor category whose stable category $\bK=\underline{\bT}$ satisfies (wfg) and central cohomological support is a weak support datum, then 
$\rho$ takes values in $\Proj C^\bullet_\bK$ and $\eta : \Proj C^\bullet_\bK \to \Spc \bK$ is a right inverse of $\rho : \Spc \bK \to \Proj C^\bullet_\bK$.}}
\item[(c)] {\em{If $\bK$ is an M$\Delta$C for which the Tate central cohomological support is a weak support datum and $\bK$ satisfies the (wTfg) condition, 
then $\widehat \eta  : \Spech \widehat C^\bullet_\bK \to \Spc \bK$ is a right inverse of $\widehat \rho  : \Spc \bK \to \Spech \widehat C^\bullet_\bK$.}}
\end{enumerate}
The distinction between parts (a) and (b) is that in the more general setting of part (a) it is not claimed that the 
map $\rho : \Spc \bK \to \Spech C^\bullet_\bK$ avoids the irrelevant ideal of $C^\bullet_\bK$, while in the 
important setting of part (b) it does, which in light of part (a), leads to the invertibility result. As always, this seemingly simple fact is the hardest 
to prove.
\subsection{\ } 
\label{1.7}
Denote by $\ThickId(\bK)$ the set of all thick two-sided ideals of an M$\Delta$C. It is an ordered monoid 
with the operation $\bI, \bJ \mapsto \langle \bI \otimes \bJ \rangle$ and the inclusion partial order. The set 
$\mc{X}_{sp}(X)$ of specialization closed subsets of a topological space $X$ is also an ordered 
monoid with the operation of intersection and the inclusion partial order. 

Our second auxiliary theorem shows that the surjective continuous maps from the main Theorem B 
are homeomorphims and classifies thick tensor ideals, under stronger assumptions than those in Theorem C.
\medskip

\noindent
{\bf{Theorem D.}} {\em{Let $\bK$ be an M$\Delta$C, which is the compact part of a compactly generated  M$\Delta$C, $\widetilde{\bK}$.

If $\bK$ satisfies the (wfg) condition, $\Proj C^\bullet_\bK$ is a Zariski space and 
the central cohomological support of $\bK$ has an extension to a faithful extended weak support datum $\widetilde{\bK} \to \mc{X}(\Proj C^\bullet_\bK)$
(see Section \ref{2.3} for definitions), then 
the following hold:}}
\begin{enumerate}
\item[(a)] {\em{The maps $\rho$ and $\eta$ are inverse homeomorphisms
\[
\Spc \bK \mathrel{\mathop{\rightleftarrows}^{\rho}_\eta} \Proj C^\bullet_\bK.
\]
}}
\item[(b)] {\em{The map 
\[
\Phi_{W_C} : \ThickId(\bK) \to \mc{X}_{sp}(\Proj C^\bullet_\bK)
\]
is an isomorphism of ordered monoids.}}
\end{enumerate}

{\em{If $\bK$ satisfies the (wTfg) condition, $\Spech \widehat C^\bullet_\bK$ is a Zariski space and 
the central cohomological support of $\bK$ has an extension to a faithful extended weak support datum $\widetilde{\bK} \to \mc{X}(\Spech \widehat C^\bullet_\bK)$,
then the following hold:}}
\begin{enumerate}
\item[(c)] {\em{The maps $\widehat{\eta}$ and $\widehat{\rho}$ are inverse homeomorphisms
\[
\Spc \bK \mathrel{\mathop{\rightleftarrows}^{\widehat \rho}_{\widehat \eta}} \Spech \widehat C^\bullet_\bK.
\]
}}
\item[(d)] {\em{The map 
\[
\Phi_{\widehat W_C} : \ThickId(\bK) \to \mc{X}_{sp}(\Spech \widehat C^\bullet_{\bK})
\]
is an isomorphism of ordered monoids.}}
\end{enumerate}

By the discussion in Section \ref{1.5},
the validity of the Etingof--Ostrik conjecture implies that the conclusions in Theorem D(a)-(b) holds for the stable category of each finite tensor category $\underline{\bT}$
for which $\Proj C^\bullet_{\underline{\bT}}$ is a Zariski space and the central cohomological support 
has an extension to a faithful extended weak support datum $\underline{\Ind{\bT}} \to \mc{X}(\Proj C^\bullet_{\underline{\bT}})$, 
see Appendix \ref{App} on background on the stable category of the indization $\underline{\Ind{\bT}}$ of $\bT$.

Theorem D recovers the well-known classifications of thick ideals and Balmer spectra in the symmetric and braided cases of stable categories of finite group schemes \cite{FP1, Balmer1}, Lie superalgebras \cite{BKN}, as well as perfect derived categories of topologically Noetherian schemes \cite{Thomason1,Balmer1}. 

We should mention connections of our work with \cite{NP} (via correpondence with Cris Negron). 
Let ${Z}^{\bullet}_{\bT}$ be the image of the restriction map from the cohomology of the Drinfeld center to the cohomology for ${\bT}$ which 
was under consideration in \cite{NP}. In particular, it was pointed out to that one expects 
$\Proj C^\bullet_{\underline{\bT}}$ to be a more precise model for $\Spc \bT$ when compared to $\Proj {Z}^{\bullet}_{\bT}$. 
Moreover, the algebra $C^\bullet_{\underline{\bT}}$ should be easier to compute as opposed to  ${Z}^{\bullet}_{\bT}$
(cf. the results in Section~\ref{sec:examples} vs. the discussions in \cite[Section 10.5]{NP}).

\medskip
\noindent
{\bf{Examples.}} Theorems B and C have relatively mild assumptions that are not hard to verify in wide generality. 
In Section \ref{sec:ex-ThmD} we demonstrate that even the most restrictive assumptions of Theorem D can be 
verified in a uniform fashion for wide classes of Hopf algebras using our new concept of the categorical center of a cohomology ring
and the results about it, obtained in the first part of the paper.
{\em{We prove that, if the assumptions in Theorem D are satisfied for a finite dimensional Hopf algebra, 
then they are satisfied for the Plavnik--Witherspoon co-smash product \cite{PW} of this Hopf algebra and 
the coordinate ring of a finite groups.}} Based on this theorem, we prove that the assumptions in Theorem D are satisfied for the 
Benson--Witherspoon Hopf algebras, and the coordinate rings of all finite group schemes.
This conveniently recovers the previous classification results for wide classes of finite tensor categories 
that were treated on a case-by-case basis with various methods. 
\subsection{\ }
\label{1.8}
We conjecture that the conclusions of parts (a)-(b) of Theorem D hold for the stable categories of all FTCs.
\medskip

\noindent
{\bf{Conjecture E.}} {\em{For every finite tensor category $\bT$, the following hold:}}
\begin{enumerate}
\item[(a)] {\em{The continuous map 
\[
\rho : \Spc \underline \bT \to \Proj C^\bullet_{ \underline \bT}
\]
is a homeomorphism.
}}
\item[(b)] {\em{The monoids $\ThickId(\underline \bT)$ and $\mc{X}_{sp}(\Proj C^\bullet_{\underline \bT})$ are isomorphic.}}
\end{enumerate}
\medskip

It is highly anticipated that the conditions in Theorem D to validate Conjecture E will hold in general. In fact, the conditions in Theorem D that an M$\Delta$C, $\bK$ satisfies 
the (wfg) condition and $\Proj C^\bullet_\bK$ is a Zariski space are natural assumptions that have appeared in the literature. In Appendix B we show that the condition that
the central cohomological support of $\bK$ has an extension to a faithful extended weak support datum $\widetilde{\bK} \to \mc{X}(\Proj C^\bullet_\bK)$
is also natural. Namely, we show that the Balmer 
support of the compact part $\bK$ of a compactly generated M$\Delta$C, $\widetilde{\bK}$,
always has an extension to an extended weak support datum. Note that the conclusion of parts (a) and (b) of Theorem D means 
that the former and latter supports coincide on thick ideals of $\bK$, which would imply that the existence of a faithful extension to $\widetilde{\bK}$ of 
the former support implies the existence of a faithful extension to $\widetilde{\bK}$ of the latter support.

\subsection{Executive Summary} 
\label{1.9}
The key conditions entering in Theorems B-D are as follows
\begin{itemize}
\item[(i)] The validity of the Etingof--Ostrik conjecture;
\item[(ii)] The weak finite generation condition; 
\item[(ii')] The weak Tate finite generation condition; 
\item[(iii)] The central cohomological support $W_C : \bK \to \Proj C_{\bK}^\bullet$ is a weak support datum map; 
\item[(iii')] The Tate central cohomological support $\widehat W_C : \bK \to \Spech \widehat C_{\bK}^\bullet$ is a weak support datum map.
\end{itemize} 
The implications of Theorems B-D are as follows:

{\em{Stable categories of finite tensor categories $\underline{\bT}$:}} for such categories the conclusion in Theorem C(a) is valid under the Etingof--Ostrik conjecture 
and that in Theorem D(a) is valid under the Etingof--Ostrik conjecture and condition (iii). We expect that condition (iii) holds in broad generality but 
it is open in key situations, such as the finite dimensional module categories of small quantum groups.
Theorem B(a) holds without any assumptions and can be helpful in studying the structure of the categorical center $C_{\underline{\bT}}^\bullet$ 
of the cohomology ring $R_{\underline{\bT}}^\bullet$, starting from the topology of $\Spc \underline{\bT}$. For instance,
prove noetherianity of $\Spc \underline{\bT}$ and then extend it to $C_{\underline{\bT}}^\bullet$, as a step to approaching 
the Etingof--Ostrik conjecture for arbitrary finite tensor categories $\bT$.

{\em{The Tate side}}: the conditions (ii') and (iii') entering in Theorems B(c), C(b) and D(b) are strong conditions on $\underline{\bT}$. The question of what finite tensor categories 
satisfy them is currently under intense study \cite{BC,CCM,N} and even the case of module categories of finite groups is still open.
In these situations Theorems B(c), C(b) and D(b) give important information about cohomological support of elements of $\bT$ computed with 
respect to Tate cohomology. 

{\em{The classification of thick tensor ideals:}} Theorems B-D with Conditions (i)-(iii') not only lead to a resolution of Problem A but also provide an explicit description of the 
Balmer spectrum of the stable categories of finite tensor categories and their thick ideals, which is of considerable interest on its own. 

{\em{Monoidal triangulated categories $\bK$ which are the compact parts of M$\Delta$Cs admitting arbitrary set indexed coproducts:}}
In this general situation, the conclusion in Theorem B(a) is valid under no assumptions, while Theorems B(b), C(a) and D(a) 
rely on assumptions (ii) and (iii). We expect that these conditions are satisfied in far greater generality than stable categories of finite tensor categories, 
thus yielding strong results towards the resolution of Problem A in full generality.
On the Tate side, Theorems B(c), C(b) and D(b) rely on assumptions (ii') and (iii'), which are more restrictive assumptions as pointed out even 
in the case of stable categories. This direction has attracted a lot of attention recently \cite{BC,C,CCM,Ng}. 
\medskip
\\
\noindent
{\bf Acknowledgements.} We acknowledge Srikanth Iyengar, Peter J{\o}rgensen, Henning Krause, and Jon Kujawa for helpful discussions and correspondence. The authors also thank Greg Stevenson, Cris Negron, and Paul Balmer for their comments on an earlier version of this manuscript. 

\section{Background on noncommutative Balmer spectra}
\label{noncommBalmerspec}
Recall that a monoidal triangulated category (M$\Delta$C for short) is a triangulated category $\bK$ with a biexact monoidal structure such that $\bK$ is a suspended monoidal category, in the language of \cite{SA}. We will assume throughout that the endomorphism ring of the unit object of $\bK$ is isomorphic to the base field $k$. 
In this section we collect some background material on the noncommutative Balmer spectrum of an M$\Delta$C that will be used in the paper, 
see \cite{Balmer1,BKS,NVY1,NVY2}. 
\subsection{The  Balmer spectrum of an M$\Delta$C}
A {\em{(two-sided) thick ideal}} $\bI$ of an M$\Delta$C, $\bK$, is a full triangulated subcategory closed under direct summands such that $\bI$ satisfies the ideal condition: for each $A \in \bI$ 
and $B \in \bK$, the objects $A \otimes B$ and $B \otimes A$ are both in $\bI$. If only $A \otimes B$ (respectively $B \otimes A$) is required to be in $\bI$, then $\bI$ is a {\em{right}} 
(respectively {\em{left}}) {\em{thick ideal}} of $\bK$. Given a collection of objects $\mc{S}$, the thick two-sided ideal generated by $\mc{S}$ will be denoted $\langle \mc{S} \rangle$. 

A two-sided thick ideal $\bP$ of $\bK$ is called {\em{prime}} if $\bI \otimes \bJ \subseteq \bP$ implies that one of $\bI$ or $\bJ$ is contained in $\bP$, 
for all thick ideals $\bI$ and $\bJ$. If $\bP$ satisfies the stronger condition that $A \otimes B \in \bP$ implies $A$ or $B$ is in $\bP$ for all objects $A$ 
and $B$ of $\bK$, then we call $\bP$ {\em{completely prime}}. 
The collection of all prime ideals of $\bK$ is the Balmer spectrum of $\bK$, denoted $\Spc \bK$, as a topological space under the Zariski topology where closed sets are defined 
to be arbitrary intersections of the base of closed sets
$$V (A) = \{ \bP \in \Spc (\bK) : A \not \in \bP \}.$$ 
The map $V$ defined above, which sends objects of $\bK$ to closed sets of $\Spc \bK$, will be referred to as the {\em{Balmer support}}. 

Recall that a subset $\mc{M}$ of $\bK$ is called {\em multiplicative} if all objects of $\mc{M}$ are nonzero, and if $A$ and $B$ are in $\mc{M}$ then so is $A \otimes B$. 
The Balmer spectrum of an M$\Delta$C, $\bK$, is always nonempty by the following result (\cite[Theorem 3.2.3]{NVY1}):

\bth{maximal} Suppose $\mc{M}$ is a multiplicative subset of a M$\Delta$C, $\bK$, and suppose $\bI$ is a proper thick ideal of $\bK$ which intersects $\mc{M}$ trivially. 
If $\bP$ is maximal element of the set 
\[
X(\mc{M}, \bI):=\{\bJ \textrm{ a thick ideal of }\bK: \bJ \supseteq \bI, \bJ \cap \mc{M}= \varnothing \},
\]
then $\bP$ is a prime ideal of $\bK$. \eth

Using Zorn's lemma and the fact that $X(\mc{M}, \bI)$ is always nonempty for any $\mc{M}$ and $\bI$ as in the hypothesis of \thref{maximal}, we conclude given such a multiplicative subset $\mc{M}$ and thick ideal $\bI$, there always exists at least one prime ideal $\bP$ containing $\bI$ and disjoint from $\mc{M}$. 
\subsection{Universality of the Balmer support}
Recall the definition of a weak support datum from Section \ref{1.6}. The Balmer support  
\[
V : \bK \to \mc{X}_{cl} (\Spc \bK)
\]
is a weak support datum, \cite[Lemma 4.1.2]{NVY1}. In \cite[Theorem 4.5.1]{NVY1} it was proved that it is a final object in the category of all 
weak support data for $\bK$:
\bth{weak-univ}
Let $\bK$ be an M$\Delta$C and $\sigma$ be a weak support datum 
from objects of $\bK$ to subsets of a topological space $X$, such that $\Phi_{\sigma}(\langle A \rangle)$ is closed for any object $A$ of $\bK$. 
Then there is a unique map $\eta_{\sigma}: X \to \Spc \bK$ satisfying $\Phi_{\sigma} (\langle A \rangle ) = \eta_{\sigma}^{-1} (V(A))$, defined by
$$\eta_\sigma(x) = \{ A \in \bK : x \not \in \Phi_{\sigma} (\langle A \rangle) \}$$ for all $x \in X$. The map $\eta_\sigma$ is continuous.
\eth
\subsection{Reconstruction of the noncommutative Balmer spectrum}
\label{2.3}
An M$\Delta$C, $\bK$, is said to be {\em{compactly generated}} if it is closed under arbitrary set indexed coproducts, 
the tensor product preserves set indexed coproducts, $\bK$ is compactly generated as a triangulated category, 
the tensor product of compact objects is compact, $\unit$ is a compact object, and every compact object is rigid (cf. \cite[Definition 2.10.11]{EGNO}). 
The full subcategory of compact objects of $\bK$, denoted by $\bK^c$, is an M$\Delta$C; we say that $\bK^c$ is the compact part of $\bK$.

Recall that for a topological space $X$, denote by $\mc{X}(X)$ the collection of all of its subsets.
An {\em{extended weak support datum}} for a compactly generated M$\Delta$C, $\bK$, is a map $\widetilde \sigma : \bK \to \mc{X}(X)$ 
that satisfies properties (1), (3)-(4) in Section \ref{1.6} and the following properties (replacing conditions (2) and (5), respectively):
\smallskip

\begin{enumerate}
\item[(2')] $\widetilde \sigma(\bigoplus_{i\in I} A_i)=\bigcup_{i \in I} \widetilde \sigma(A_i )$, for all  $A_i \in \bK$;
\item[(5')] $\bigcup_{D \in \bK^c} \widetilde \sigma(A \otimes D \otimes C) = \widetilde \sigma(A) \cap \widetilde \sigma(C)$ for all $A \in \bK, C \in \bK^c$.
\end{enumerate}
\smallskip

This is a slightly weaker assumption that the notion of extended weak support datum used in \cite{NVY1,NVY2}.
We say that such a support map satisfies the {\em{faithfulness property}} if
\[
\Phi_{\widetilde \sigma}(\langle M \rangle_{\bK})=\varnothing \Leftrightarrow M=0, \quad \forall M \in \bK
\] 
and the {\em{realization property}} if for every $W \in {\mathcal X}_{cl}(X)$ there exists  $M\in \bK^{c}$ such that $\Phi_{\widetilde \sigma}(\langle M\rangle_{\bK^c})=W$. 
Here, $\langle M\rangle_{\bK^c}$ and $\langle M\rangle_{\bK}$ denote the thick ideals of $\bK^c$ and $\bK$ generated by an object $M$ in one of them.
We say that $\widetilde \sigma : \bK \to \mc{X}(X)$ is an {\em{extension}} of $\sigma : \bK^c \to \mc{X}_{sp}(X)$ if 
\[
\Phi_{\widetilde \sigma}(\langle M \rangle_{\bK^c}) = \Phi_{\sigma}(\langle M \rangle_{\bK^c}) \quad \mbox{for all} \quad 
M \in \bK^c. 
\]


We will need the following reconstruction result for the noncommutative Balmer spectrum, which is a slightly more general version of \cite[Theorem 6.2.1]{NVY1}.
Its proof is identical to that of \cite[Theorem 6.2.1]{NVY1}.

\bth{reconstr}
Assume that $\bK$ is a compactly generated M$\Delta$C and $\sigma : \bK \to \XX$ is an extended weak support datum
for a Zariski space $X$ such that $\Phi_\sigma (\langle C \rangle)$ is closed for every compact object $C$.
Suppose in addition that $\sigma$ satisfies the faithfulness and realization properties. Then the following hold:
\begin{enumerate}
\item[(a)] The map $\eta_\sigma: X \to \Spc \bK^c$ is a homeomorphism.
\item[(b)] The map $\Phi_\sigma$ is an isomorphism of ordered monoids between the set of thick ideals of $\bK^c$, equipped with the operation $\bI, \bJ \mapsto \langle \bI \otimes \bJ \rangle$ 
and the inclusion partial order, and the set $\mc{X}_{sp}(X)$ of specialization closed subsets of $X$, equipped with the operation of intersection and the inclusion partial order.
\end{enumerate}
\eth


\section{Nine algebras}
\label{sec:five-algebras}
In this section we describe the algebras of homological origin that 
will play a role in this paper and exactly where the categorical centers fit into this picture.
We compare the latter rings to various families of algebras that have appeared in the literature.
\subsection{Relation between the categorical center and the graded center of an M$\Delta$C}
Let $\bK$ be an M$\Delta$C with a fixed set of generators $\mc{K}$. 
Recall from the introduction that the categorical center $\widehat C^\bullet_{\bK}$ of the 
Tate cohomology ring $\widehat R^\bullet_{\bK}$ is the graded subring spanned by all $g \in \Hom_{\bK}(\unit, \Sigma^i \unit)$
such that the diagram \eqref{diag-M-g} commutes for all objects $M\in {\mathcal K}$. 

A pathological problem in this situation is that, if $M_1 \to M_2 \to M_3 \to \Sigma M_1$ is a distinguished triangle and the diagram 
\eqref{diag-M-g} commutes for $M = M_1$ and $M=M_3$, then it does not necessarily commute for $M = M_2$. 
This comes from the non-uniqueness of the extension morphism in the TR4 axiom for triangulated categories.
Consequently, the 
algebra $\widehat C^\bullet_{\bK}$ depends on the choice of the collection $\mc{K}$ and is bigger than 
the algebra that is obtained in this way when we choose $\mc{K}$ to be the set of all objects of $\bK$. Because of this, 
we choose $\mc{K}$ to be as small as possible, and in particular, in the case of the
stable category $\underline{\bT}$ of a finite tensor category $\bT$, 
${\mathcal K}$ will be always chosen to be the (finite) set of simple 
objects of $\bT$. This is illustrated in detail in Section \ref{sec:examples}.

The {\em{graded center}} of $\bf{K}$, as studied in
\cite{BIK1,BF}, is the (graded commutative) algebra $Z^\bullet({\bK})$ whose degree $n$ component 
consists of natural transformations 
\[
\eta : \id_{\bK} \to \Sigma^n \quad \mbox{such that} \quad \eta \Sigma = (-1)^n \Sigma \eta.
\]  
By \cite[Proposition 2.1]{BKSS}, there are two injective homomorphisms 
$\mc{L}, \mc{R} : \widehat R^\bullet_{\bK} \hookrightarrow Z^\bullet({\bK})$, 
which send $g \in \Hom_{\bf{K}} ({\bf{1}}, \Sigma^n {\bf{1}})$ 
to
\[
\begin{tikzcd}
M \arrow[r, "\cong"] & {\bf{1}} \otimes M \arrow[r, "g \otimes \id_M"] & \Sigma^n {\bf{1}} \otimes M \arrow[r, "\cong"] & \Sigma^n M  \quad \mbox{and}
\\
M \arrow[r, "\cong"] & M \otimes {\bf{1}} \arrow[r, "\id_M \otimes g"] & M \otimes \Sigma^n {\bf{1}} \arrow[r, "\cong"] & \Sigma^n M,   
\end{tikzcd}
\]
respectively. By way of definition, the categorical center of the Tate cohomology ring $\widehat R^\bullet_{\bK}$  is given by
\[
\widehat C^\bullet_{\bK} := \{ g \in \widehat R^\bullet_{\bK} \mid \mc{L}(g) = \mc{R}(g) \; \; \mbox{on all objects} \; \; M \in \mc{K} \}.
\]
This realizes $\widehat C^\bullet_{\bK}$ as a subalgebra of the graded center $Z^\bullet({\bK})$, 
which is generally much smaller than $Z^\bullet({\bK})$ since the Tate cohomology ring 
$\widehat R^\bullet_{\bK}$ is itself embedded in $Z^\bullet({\bK})$.


\subsection{The Drinfeld center}
Recall that the {\em{Drinfeld center}} of a monoidal category $\bT,$ denoted $\bZ(\bT)$, is the category with objects given by pairs $(A, \gamma)$, where $A$ is an object of $\bT$ and 
$\gamma$ is a half-braiding, that is, a natural isomorphism 
\[
\gamma_M : A \otimes M \cong M \otimes A, \quad M \in \bT
\]
satisfying the associativity condition \cite[Eq.\ (7.41)]{EGNO}. It is again a monoidal category with unit object
$(\unit, \beta)$, where $\unit$ is the unit object of $\bT$ and $\beta$ is the structure isomorphism $\beta_M: \unit \otimes M \cong M \cong M \otimes \unit$ in $\bT$. 
A morphism $(A, \gamma) \to (A', \gamma')$ in $\bZ(\bT)$ is defined to be a morphism $g: A \to A'$ in $\bT$ such that the diagram 
\begin{center}
\begin{tikzcd}
A \otimes M \arrow[r, "g \otimes \id_M"] \arrow[d, "\gamma_M"] & A' \otimes M \arrow[d, "\gamma'_M"] \\
M \otimes A \arrow[r, "\id_M \otimes g"]                       & M \otimes A'                       
\end{tikzcd}
\end{center}
commutes for all objects $M$ of $\bT$, see \cite[Sect.\ 7.13-14]{EGNO} for details. 
For example, if $\bT= \modd (H)$ for a finite dimensional Hopf algebra $H$, then $\bZ(\bT)$ is equivalent to the category $\modd(D(H))$, where $D(H)$ is the Drinfeld double of $H$ 
(cf. \cite[Proposition 7.14.6]{EGNO}). 

\subsection{Cohomology rings of M$\Delta$Cs and their categorical centers}
The Drinfeld center $\bZ(\bK)$ of an M$\Delta$C, $\bK$,
is a monoidal additive category equipped with the autoequivalence $(A, \gamma) \mapsto (\Sigma A, \gamma')$, where $\gamma'$ is the natural isomorphism defined by the composition of $\Sigma \gamma$ with the M$\Delta$C structure isomorphisms:
\begin{center}
\begin{tikzcd}
(\Sigma A) \otimes M \arrow[d, "\cong"'] \arrow[r, "\gamma'_M"] & M \otimes (\Sigma A)                    \\
\Sigma(A \otimes M) \arrow[r, "\Sigma(\gamma_M)"']              & \Sigma(M \otimes A) \arrow[u, "\cong"']
\end{tikzcd}
\end{center} 
(We are not aware of an argument showing that $\bZ(\bK)$ is necessarily a triangulated category.) 

Thus, we have another canonical algebra attached to $\bK$: 
{\em{the Tate cohomology ring of the Drinfeld center}} $\bZ(\bK)$, 
\[
\widehat R^\bullet_{\bZ(\bK)}:= \bigoplus_{i \in \mathbb{Z}} \Hom_{\bZ(\bK)}((\unit, \beta), \Sigma^i (\unit, \beta)).
\]
It is easy to verify that forgetful functor $\bZ(\bK) \to \bK$ induces an injective homomorphism
\begin{equation}
\label{eq:inj-D-to-C}
\widehat R^\bullet_{\bZ(\bK)} \hookrightarrow \widehat C^\bullet_{\bK} \quad \mbox{given by} \quad 
g \in \Hom_{\bZ(\bK)}((\unit, \beta), \Sigma^i (\unit, \beta)) \mapsto g \in \Hom_{\bK}(\unit, \Sigma^i \unit),
\end{equation}
which is an isomorphism if the collection $\mc{K}$ is chosen to be the set of all objects of $\bK$.

Similarly, we have the {\em{cohomology ring of the Drinfeld center}} $\bZ(\bK)$, 
\[
R^\bullet_{\bZ(\bK)}:= \bigoplus_{i \geq 0} \Hom_{\bZ(\bK)}((\unit, \beta), \Sigma^i (\unit, \beta)),
\]
and the injective homomorphism \eqref{eq:inj-D-to-C} restricts to an injective homomorphism 
\[
R^\bullet_{\bZ(\bK)} \hookrightarrow C^\bullet_{\bK},
\]
which is an isomorphism if the collection $\mc{K}$ is chosen to be the set of all objects of $\bK$.
\subsection{Relations in the case of stable categories of finite tensor categories}
Recall that a finite tensor category $\bT$ is a tensor category which is equivalent to the category of finite 
dimensional representations of a finite dimensional $k$-algebra. See \cite[Ch. 4]{EGNO} for background on finite tensor categories.
 
Consider the important special case when $\bK$ is the stable category of a finite tensor category $\bT$ 
which is not semisimple; 
as usual the corresponding stable category will be denoted by $\underline{\bT}$. 
In this situation the collection ${\mathcal K}$ will be chosen to be the (finite) set 
of simple objects of $\mc{C}$. Furthermore, $R^\bullet_{\underline{\bT}}$ is isomorphic to the {\em{cohomology ring}} of $\bT$
\[
R_{\bT}^\bullet := \bigoplus_{i \geq 0} \Ext^i_{\bT}(\unit, \unit).
\]
This follows from the canonical isomorphism 
\[
\bigoplus_{i > 0} \Ext^i_{\bT} (A, B) \cong \bigoplus_{i > 0} \Hom_{\underline{\bT}} (A, \Sigma^i B)
\]
for all objects $A, B$ of a Frobenius category $\bT$ (see e.g., \cite[Proposition 2.6.2]{CTVZ1} for the case of group algebras, which extends to arbitrary Frobenius categories), and from the isomorphisms 
$\End_{\bT}(\unit) \cong k$ and $\End_{\underline{\bT}} (\unit) \cong k$ 
(the first isomorphism is part of the definition of a finite tensor category, the second 
one follows from the assumption that $\bT$ is not semisimple, which implies that $\unit$ is not projective).
The algebra $\widehat R^\bullet_{\underline{\bT}}$ is called the {\em{Tate cohomology ring}} of $\bT$. 

In the special case of the stable module category $\underline{\modd}(H)$ of a finite dimensional Hopf algebra $H$, 
$R^\bullet_{\underline{\modd}(H)}$ and $\widehat R^\bullet_{\underline{\modd}(H)}$ are isomorphic to 
the cohomology ring and the Tate cohomology ring of $H$, respectively.

The Drinfeld center $\bZ(\bT)$ of a finite tensor category $\bT$ is also a finite tensor category \cite[Proposition 7.13.8]{EGNO}, 
and hence we can additionally form its stable category $\underline{\bZ(\bT)}$. 

The categorical centers $C^\bullet_{\underline{\bT}} \subseteq \widehat C^\bullet_{\underline{\bT}} $ 
are closely related to two additional algebras of interest: {\em{the cohomology ring of the Drinfeld center}} $\bZ(\bT)$, 
\[
R_{\bZ(\bT)}^\bullet := \bigoplus_{i\geq 0} \Ext^i_{\bZ(\bT)} ((\unit, \beta), (\unit, \beta)) \cong R_{\underline{\bZ(\bT)}}^\bullet
\]
{\em{and the Tate cohomology ring of the stable category of the Drinfeld center of $\bT$}}, that is, $\underline{\bZ(\bT)}$: 
\[
\widehat R_{\underline{\bZ(\bT)}}^\bullet := \bigoplus_{i \in \mathbb{Z}} \Hom_{\underline{\bZ(\bT)}} ( (\unit, \beta), \Sigma^i (\unit, \beta)).
\]

By \cite[Proposition 2.1.1]{V}, there is a monoidal triangulated functor $\overline{F}: \underline{\bZ(\bT)} \to \underline{\bT}$ which extends the usual forgetful functor $F: \bZ(\bT) \to \bT$. If $(A, \gamma) \in \underline{\bZ(\bT)}$-- i.e., $\gamma: A \otimes - \cong - \otimes A$ is a half-braiding for $\bT$-- then it is straightforward to verify that $\gamma$ also defines a half-braiding on $\underline{\bT}$. If $g: A \to A'$ defines a morphism $(A, \gamma) \to (A', \gamma')$ in $\underline{\bZ(\bT)}$, then the diagram 
\begin{center}
\begin{tikzcd}
A \otimes M \arrow[r, "g \otimes \id_M"] \arrow[d, "\gamma_M"] & A' \otimes M \arrow[d, "\gamma'_M"] \\
M \otimes A \arrow[r, "\id_M \otimes g"]                       & M \otimes A'                       
\end{tikzcd}
\end{center}
commutes in $\bT$ for any object $M$, and hence commutes in $\underline{\bT}$; thus $g$ also defines a morphism $\check{F}(A,\gamma) \to \check{F}(A', \gamma')$. Therefore, $\overline{F}: \underline{\bZ(\bT)} \to \underline{\bT}$ factors as
\begin{center}
\begin{tikzcd}
\underline{\bZ(\bT)} \arrow[r, "\check{F}"'] \arrow[rr, "\overline{F}", bend left] & \bZ(\underline{\bT}) \arrow[r] & \underline{\bT},
\end{tikzcd}
\end{center}
where the functor $\bZ(\underline{\bT}) \to \underline{\bT}$ is the forgetful functor. The functor $\check{F}$ is additive and monoidal, and commutes with the respective shift functors on $\bZ(\underline{\bT})$ and $\underline{\bZ(\bT)}$, and thus induces a homomorphism of graded rings $$i: \widehat R^\bullet_{\bZ(\bT)} \cong \widehat R^\bullet_{\underline{\bZ(\bT)}} \to  \widehat R^\bullet_{\bZ(\underline{\bT})},$$ 
sending $g: (\unit, \beta) \to \Sigma^i (\unit, \beta)$ to $\check{F}(g): (\unit, \beta) \to \Sigma^i (\unit,\beta)$.


\ble{im-i-abull}
The image of the composition of ring homomorphisms $\widehat R^\bullet_{\underline{\bZ(\bT)}} \to \widehat R^\bullet_{\bZ(\underline{\bT})} \to \widehat R^\bullet_{\underline{\bT}}$ is contained in $\widehat C^\bullet_{\underline{\bT}}$.
\ele

\begin{proof}

Let $g \in \widehat R^\bullet_{\underline{\bZ(\bT)}}$ be a map in $\underline{\bZ(\bT)}$ denoted by $g:(\unit, \beta) \to \Sigma^i (\unit, \beta).$ As a morphism in the stable category of $\bZ(\bT)$, $g$ corresponds to an equivalence class of morphisms in $\bZ(\bT)$ from $(\unit, \beta) \to \Sigma^i(\unit, \beta)$. Since $g$ corresponds to a morphism in the Drinfeld center of $\bT$, by definition the diagram
\begin{center}
\begin{tikzcd}
\unit \otimes A \arrow[r, "g \otimes \id_A"] \arrow[d] & \Sigma^i \unit \otimes A \arrow[d] \\
A \otimes \unit \arrow[r, "\id_A \otimes g"]           & A \otimes \Sigma^i \unit          
\end{tikzcd}
\end{center}
commutes, for all $A \in \bT$, and where the vertical maps are the half-braidings associated to $\unit$ and $\Sigma^i \unit$, respectively. Since the half-braidings associated to these objects are the same as the structure maps in the monoidal triangulated structure, the commutativity of this diagram implies that the image of $g$ under this composition, in other words $\overline{F}(g)$, is in $\widehat C^\bullet_{\underline{\bT}}$.

 \end{proof}

For $\underline{\bT}$ the stable module category of any finite tensor category $\bT$, we have now exhibited a collection of nine algebras (up to isomorphism):
\begin{equation}
\begin{tikzcd}
R^\bullet_{\bZ(\bT)}\cong R^\bullet_{\underline{\bZ(\bT)}} \arrow[r, "i"] \arrow[rd, hook] & R_{\bZ(\underline{\bT})}^\bullet \arrow[rd, hook] \arrow[r,hook] & C^\bullet_{\underline{\bT}} \arrow[rd, hook] \arrow[r, hook] & R^\bullet_{\bT} \cong R^\bullet_{\underline{\bT}} \arrow[rd, hook] &                                                                          &                            \\
                                                                                           & \widehat R^\bullet_{\underline{\bZ(\bT)}} \arrow[r, "i"]    & \widehat R_{\bZ(\underline{\bT})}^\bullet \arrow[r,hook]        & \widehat C^\bullet_{\underline{\bT}} \arrow[r, hook]                 & \widehat R^\bullet_{\underline{\bT}} \arrow[r, "{\mc{R}, \mc{L}}", hook] & Z^\bullet(\underline{\bT})
\end{tikzcd}
\label{eq:nine-algs}
\end{equation}
\bre{skipK} From now on, when the underlying M$\Delta$C ($\bK$ or $\underline{\bT}$) is clear from the discussion, we will omit the subscript in the notation for the above algebras and will simply write
\[
C^\bullet, \widehat C^\bullet, R^\bullet, \widehat R^\bullet.
\] 
Analogously, we will write $\End^\bullet(A)$, $\hEnd^\bullet(A)$,  $\Hom^\bullet(A,B)$, $\hHom^\bullet(A, B)$ without subscripts $\bK$ or $\underline{\bT}$.
\ere 

\section{Categorical centers of cohomology rings and fixed point subrings}
\label{sec:fixed-pt}
In this section we show that, in certain important situations, the categorical centers $C_{\bK}^\bullet$ and $\widehat C_{\bK}^\bullet$ of the 
cohomology rings $R_{\bK}^\bullet$ and $\widehat R_{\bK}^\bullet$ can be interpreted as fixed point subrings with respect to a natural 
group action.

\subsection{The group of endotrivial objects of an M$\Delta$C}
Let $\bK$ be an arbitrary monoidal category. Recall that an object $M \in \bK$ is called {\em{left dualizable}} if there exists an object $M^* \in \bK$ 
and evaluation, coevaluation maps
\[
\ev: M^* \otimes M \to \unit, \quad
\coev: \unit \to M \otimes M^*,
\] 
such that the compositions
\[
M \xrightarrow{\coev \otimes \id} M \otimes M^* \otimes M \xrightarrow{ \id \otimes \ev} M \quad \mbox{and} \quad 
M^* \xrightarrow{\id \otimes \coev} M^* \otimes M \otimes M^* \xrightarrow{\ev \otimes \id} M^*
\]
are the identity maps on $M$ and $M^*$, respectively. The left dual object $M^*$ is unique up to a unique isomorphism, 
\cite[Proposition 2.10.5]{EGNO}. In a similar way one defines the notions of a {\em{right dualizable}} object $M$ and its {\em{right dual}} ${}^*M$, 
see \cite[Definition 2.10.2]{EGNO}. An object $M \in \bK$ is called {\em{rigid}} if it is both left and right dualizable.

An object $M \in \bK$ is called {\em{invertible}} if it is left dualizable and the maps $\ev$ and $\coev$ are isomorphisms, cf. 
\cite[Definition 2.11.1]{EGNO}. 
\ble{invert-obj} Every invertible object $M$ of a monoidal category $\bK$ is rigid with ${}^* M \cong M^*$. The set of invertible
objects of $\bK$ forms a group under $\otimes$, with the inverse of $M \in \bK$ equal to $M^*$. 
\ele
\begin{proof} A result of this type was proved in the case of rigid monoidal categories in \cite[Proposition 2.11.3]{EGNO}. 
We provide brief details because of the more general setting that is considered here. For an invertible 
object $M \in \bK$, we have the maps 
\[
\coev^{-1}: M \otimes M^* \to \unit, \quad
\ev^{-1}: \unit \to M^* \otimes M, 
\] 
and one easily checks that $M^*$ and this pair of maps satisfy the axioms for right dual object of $M$; thus $M^*$ is also a right dual of $M$.  
This, in particular implies that $M^*$ is also invertible with left/right dual given by $M$ and that $M^{**} \cong M$. The rest of the lemma is 
now straightforward. 
\end{proof}

With the previous lemma, one can define the group of endotrivial modules for any rigid monoidal category. As far as the authors know this definition has not 
appeared in the context of non-symmetric monoidal categories. 

\bde{endo}Let $\bK$ be a monoidal category where every object is rigid. For any object in $M\in \bK$, let $[M]$ note the isoclass of objects defined by $M$. Let 
$T(\bK)$ be the group consisting of the classes $[M]$ for invertible objects $M$ with the product structure $[M]\cdot [N]=[M\otimes N]$,
identity object $[\unit]$, and inverses $[M]^{-1} := [M^*]$.
We will call $T(\bK)$ the {\em group of endotrivial objects} for $\bK$,
extending the standard terminology from representations of finite groups. 
\ede 

For $H$ a finite dimensional Hopf algebra, one can consider the stable module category  $\underline{\modd}(H)$, and the group of endotrivial 
modules $T(H):=T(\underline{\modd}(H))$. In the case where $H$ is 
cocommutative, $T(H)$ is an abelian group. Carlson and Nakano have conjectured that $T(H)$ should be finitely generated \cite[Sec.~10, (2)]{CN}. In the important case when $H=kG$ is a group algebra of a finite group the theory of endotrivial modules has been well-studied (cf. \cite{Mazza}). Puig \cite{P} proved that $T(kG)$ is a finitely generated (abelian) group. The structure of $T(kG)$ has been computed for a various families of groups 
(cf. \cite{CMN1, CMN2, CMN3, CMN4}). When $H$ is an arbitrary Hopf algebra, an interesting question is to determine if  $T(H)$ is finitely generated, and to work out the group structure. 

\subsection{The action of the group of endotrivial modules}
For the rest of the section, let $\bK$ be an M$\Delta$C. Let $T(\bK)$ be the group of endotrivial objects of $\bK$. It acts on the Tate cohomology ring 
$\widehat R^\bullet$ by homogeneous automorphisms (and thus, on $R^\bullet$) in the following way. The action of $M \in T(\bK)$ on 
$g \in \Hom(\unit, \Sigma^i \unit)$ is the morphism $\sigma_M(g) \in \Hom(\unit, \Sigma^i \unit)$ given by the composition 
\begin{align}
&\unit \xrightarrow{\coev} M \otimes M^* \cong M \otimes \unit \otimes M^* \xrightarrow{\id_M \otimes g \otimes \id_{M^*}} M \otimes \Sigma^i \unit \otimes M^* \cong
\label{sigma-M-g}
\\
& \cong \Sigma^i( M \otimes \unit \otimes M^*) \cong \Sigma^i (M \otimes M^*) \xrightarrow{\Sigma^i \coev^{-1} }  \Sigma^i \unit.
\nn
\end{align}
\ble{comm-M-g} For an M$\Delta$C, $\bK$, and $g  \in \Hom(\unit, \Sigma^i \unit)$, the following hold:
\begin{enumerate}
\item[(a)] If the diagram \eqref{diag-M-g} commutes, then the diagram \eqref{diag-M-g} with $M$ replaced by $\Sigma^n M$ 
commutes for all $n \in \mathbb{Z}$.
\item[(b)] If the diagram \eqref{diag-M-g} commutes for $M= M_1$ and $M=M_2$, then it commutes for $M=M_1 \otimes M_2$.
\item[(c)] For an invertible $M \in \bK$, the diagram \eqref{diag-M-g} commutes if and only if $\sigma_M(g) =g$. 
\end{enumerate}
\ele
\begin{proof} Parts (a) and (b) are straightforward and are left to the reader. 

(c) If the diagram \eqref{diag-M-g} commutes and $M$ is invertible, then one can replace the second arrow in \eqref{sigma-M-g} with 
the composition of the other three arrows in the commutative diagram
\begin{center}
\begin{tikzcd}
M \otimes \unit \otimes M^* \arrow[r, "\id_M \otimes g \otimes \id_{M^*}"] \arrow[d, "\cong"] & [4em]M \otimes \Sigma^i \unit \otimes M^*  \\
\unit \otimes M \otimes M^* \arrow[r, "g \otimes \id_M \otimes \id_{M^*}"]           & \Sigma^i \unit \otimes M \otimes M^* \arrow[u, "\cong"]        
\end{tikzcd}
\end{center}
where in the second column $ \Sigma^i \unit \otimes M \otimes M^* \cong \Sigma^i(M \otimes M^*) \cong M \otimes \Sigma^i \unit \otimes M^*$.
The cancellation of $\coev$ and $\coev^{-1}$ and the bi-functoriality of the tensor product in $\bK$ gives that 
$\sigma_M(g) =g$. The opposite direction is proved in an analogous way.
\end{proof}
\leref{comm-M-g}(c) implies that for every M$\Delta$C, $\bK$, we have the inclusions of graded algebras
\[
\widehat C^\bullet \subseteq (\widehat R^\bullet)^{T(\bK) \cap [\mc{K}]} \subseteq \widehat R^\bullet \quad \mbox{and} \quad
 C^\bullet \subseteq ( R^\bullet)^{T(\bK)  \cap [\mc{K}]} \subseteq  R^\bullet,
\]
where $[\mc{K}]$ denotes the collection of isoclasses of objects of $\bK$ corresponding to the objects in $\mc{K}$. 
In the following case, the first inclusions in each of the two chains become equalities, and the categorical centers have a concrete realization as rings of fixed points. 
\bco{isotriv} If an M$\Delta$C, $\bK$, is generated by a subgroup $G$ of its group of endotrivial objects $T(\bK)$ and $[\mc{K}] :=G$, then 
\[
\widehat C^\bullet = (\widehat R^\bullet)^G \quad \mbox{and} \quad
C^\bullet = (R^\bullet)^G.
\]
\eco
\begin{proof}
The corollary follows from \leref{comm-M-g}(c).
\end{proof}
\section{Examples of categorical centers} 
\label{sec:examples}
In this section we illustrate how the categorical centers $C_{\bK}^\bullet$ and $\widehat C_{\bK}^\bullet$ of the 
cohomology rings $R_{\bK}^\bullet$ and $\widehat R_{\bK}^\bullet$ arise as fixed point subrings
for the stable module categories of important classes of Hopf algebras. For a Hopf algebra 
$A$ over a field $k$, we will denote by 
\begin{equation}
\label{Tate-coh}
\widehat \opH^\bullet(A, k) \cong \widehat R^\bullet_{\underline{\modd}(A)}
\end{equation}
its Tate cohomology ring.
\subsection{Basic Hopf algebras} A finite dimensional Hopf algebra $A$ is called {\em basic} if all its irreducible representations are 1-dimensional 
(i.e., $A$ is basic when considered as an algebra). 
In this situation, the irreducible representations are endotrivial objects of $\underline{\modd}(A)$, and the set of irreducibles is closed under tensoring. 
Furthermore, $\underline{\modd}(A)$ is generated by that subgroup of endotrivial objects. The following 
proposition gives an explicit description of the categorical centers in this case. 

\bpr{coh-cent-fix} Let $A$ be a finite dimensional basic Hopf algebra and $G$ be the set of isomorphism classes of non-zero 1-dimensional $A$-modules 
in the stable module category $\underline{\modd}(A)$ of finite dimensional $A$-modules. Then $G$ is a subgroup of $T(\underline{\modd}(A))$ 
and 
\[
\widehat C_{\underline{\modd}(A)}^\bullet = (\widehat R^\bullet_{\underline{\modd}(A)})^G \cong ( \widehat \opH^\bullet (A, k))^G
\quad \mbox{and} \quad
C^\bullet_{\underline{\modd}(A)} = (R^\bullet_{\underline{\modd}(A)})^G \cong (\opH^\bullet (A, k))^G.
\]
\epr

\begin{proof} Every non-zero 1-dimensional $A$-module is invertible in the abelian category ${\modd}(A)$ of finite dimensional $A$-modules
with inverse given by its left=right dual.
The set of those modules is closed under tensoring and taking inverses. Hence, $G$ is a subgroup of $T(\underline{\modd}(A))$. 
Since $A$ is basic, all irreducible $A$-modules are one dimensional, and thus every module in ${\modd}(A)$ has a composition 
series with one dimensional subquotients. Therefore, $G$ generates $\underline{\modd}(A)$ as a triangulated category and 
the proposition follows from \coref{isotriv}.
\end{proof}

\bex{UnipHopgalg}[{\bf Unipotent Hopf algebras}] Assume that $A$ is a finite dimensional Hopf algebra which is local when considered as an algebra; we will call such a Hopf algebra 
a {\em{unipotent Hopf algebra}}. It has a unique irreducible representation coming from its counit, and this representation is the 
unit object $\unit$ of $\modd (A)$. One can then take 
$G$ to be the image of $\unit$ in $\underline{\modd}(A)$. Since $\unit \otimes M \cong M \otimes \unit$ for all $M \in \modd (A)$, 
$G$ acts trivially on $\widehat R^\bullet_{\underline{\modd}(A)}$. \prref{coh-cent-fix}(b) implies that
\[
\widehat C_{\underline{\modd}(A)}^\bullet = \widehat R^\bullet_{\underline{\modd}(A)} \cong \widehat \opH^\bullet (A, k)
\quad \mbox{and} \quad
C^\bullet_{\underline{\modd}(A)} = R^\bullet_{\underline{\modd}(A)} \cong \opH^\bullet (A, k). 
\]
\eex

\bex{QuantBorel}{[\bf Quantum Borels]} Let ${\mathcal R}$ be an irreducible root system of rank $n$. Let $\ell$ be a positive integer and 
$\zeta$ be a primitive $\ell$th root of unity. Assume that ${\ell}$ is a positive integer such that  

(i) ${\ell}$ is odd, (ii) if ${\mathcal R}$ is of type $G_{2}$ then $3\nmid {\ell}$, and (iii) ${\ell}>h$ where $h$ is the Coxeter number for ${\mathcal R}$. 

Let $\mathfrak g$ be the complex simple Lie algebra associated with ${\mathcal R}$, and ${\mathfrak b}$ be a Borel subalgebra. For every lattice 
$\Lambda$ with ${\mathbb Z}{\mathcal R} \subseteq \Lambda \subseteq X$, where $X$ is the weight lattice, one can construct a small quantum Borel 
subalgebra $u_{\zeta,\Lambda}({\mathfrak b})$, \cite{N}.
We will follow the notation of \cite[Sect.\ 5.1]{NVY2}; the generators of $u_{\zeta,\Lambda}({\mathfrak b})$ will be denoted by $E_{\alpha_i}$ and 
$K_\mu$ for $\{\alpha_1, \ldots, \alpha_n\}$ a base of simple roots of ${\mathcal R}$ and $\mu \in \Lambda/\Lambda'$, 
where the sublattice $\Lambda'$ of $\Lambda$ is given by the display below Eq. (5.1.1) in \cite{NVY2}.

It is well known that $u_{\zeta,\Lambda}({\mathfrak b})$ is basic (see e.g., \cite[Proposition 5.3.1(a)]{NVY2}), 
so \prref{coh-cent-fix}(b) applies to it. 
By \cite[Proposition 5.3.1(c)]{NVY2}, the group $G$ in \prref{coh-cent-fix}(b) acts trivially on $R^\bullet_{\underline{\modd}(u_{\zeta,\Lambda}({\mathfrak b}))}$. Therefore, 
\[
C_{\underline{\modd}(u_{\zeta,\Lambda}({\mathfrak b}))}^\bullet = R^\bullet_{\underline{\modd}(u_{\zeta,\Lambda}({\mathfrak b}))} 
\cong \opH^\bullet(u_{\zeta,\Lambda}({\mathfrak b}), {\mathbb C}).
\]
Using the Lyndon-Hochschild-Serre spectral one shows that 
\[
\widehat R^\bullet_{\underline{\modd}(u_{\zeta,\Lambda}({\mathfrak b}))}\cong 
\widehat{\opH}^{\bullet}(u_{\zeta}({\mathfrak u}),{\mathbb C})^{u_{\zeta,\Gamma}({\mathfrak t})},
\]
where $u_{\zeta}({\mathfrak u})$ and $u_{\zeta,\Gamma}({\mathfrak t})$ are the unital subalgebras of $u_{\zeta,\Lambda}({\mathfrak b})$
generated by the elements $\{E_{\alpha_i} \mid 1 \leq i \leq r\}$ and $\{K_\mu \mid \mu \in \Lambda/ \Lambda' \}$, respectively.
Now we can invoke the argument in \cite[Proposition 5.3.1(c)]{NVY2} to show that the group $G$ acts trivially on 
$\widehat{R}^\bullet_{\underline{\modd}(u_{\zeta,\Lambda}({\mathfrak b}))}$. Therefore, 
\begin{equation}
\label{cent-Borel}
\widehat C_{\underline{\modd}(u_{\zeta,\Lambda}({\mathfrak b}))}^\bullet = \widehat R^\bullet_{\underline{\modd}(u_{\zeta,\Lambda}({\mathfrak b}))}
\cong \widehat \opH^\bullet(u_{\zeta,\Lambda}({\mathfrak b}), {\mathbb C}).
\end{equation}
\eex 

\subsection{Plavnik-Witherspoon co-smash product Hopf algebras}
\label{plav-with} We begin with a construction of families of finite dimensional Hopf algebras introduced by Plavnik and Witherspoon \cite[Section 2]{PW}. 
Let $L$ be a finite group and $A$ be a finite dimensional Hopf algebra with $L$ acting on $A$ by Hopf automorphisms. Set 
$$A_{L}=(A^{*}\# kL)^{*}.$$ 
As an algebra, $A_L \cong A \otimes k[L]$. The algebra $k[L]$ is semisimple and the central idempotents are indexed by $g\in L$. 
Therefore, the modules for $A_{L}$ are the same as $L$-graded $A$-modules. We write such modules in the form 
\begin{equation}
\label{AL-decomp}
M=\bigoplus_{x \in L} N_x \otimes \kk_x,
\end{equation}
where each $N_x$ is an $A$-module and $\kk_x$ is the 1-dimensional $k[L]$ corresponding to $x \in L$. 

There is a functor $\mc{F}: \modd(A) \to \modd(A_{L})$, defined by 
\[
N \mapsto N \otimes \kk_e
\]
on objects, where $e$ is the identity of $L$. The functor $\mc{F}$ extends to a fully faithful monoidal triangulated functor $\underline{\modd}(A) \to \underline{\modd}(A_{L})$;
this is proved similarly to \cite[Lemma 9.2.2]{NVY1}. The functor $\mc{F}$ induces an isomorphism between the Tate cohomology rings of $A_{L}$ and $A$:
\begin{equation}
\label{AL-A}
\widehat R^\bullet_{\underline{\modd}(A_{L})} \cong \widehat R^\bullet_{\underline{\modd}(A)}.
\end{equation}

Consider the collection $S$ of objects of $\modd(A_{L})$, defined as the objects of the forms 
\begin{enumerate}
\item $N \otimes \kk_e$, for $N$ a simple $A$-module, and
\item $\epsilon \otimes \kk_x$, where $\epsilon$ is the trivial $A$-module, and $x$ is an element of $L$.
\end{enumerate}
Every object of the generating set $\mc{K}$ (that is, every irreducible $A_L$-module) can be written as a tensor product of elements of $S$. The formula for the tensor product of modules of $A_{L}$, as proven in \cite{PW}, gives 
\[
(N \otimes \kk_e) \otimes (\epsilon \otimes \kk_x) \cong (N \otimes ^x \epsilon)\otimes (\kk_e \otimes \kk_x) \cong (N \otimes \epsilon) \otimes \kk_x \cong N \otimes \kk_x.
\]
This implies that the identification \eqref{AL-A} gives an embedding
\[
\widehat C^\bullet_{\underline{\modd}(A_{L})} \hookrightarrow  \widehat C^\bullet_{\underline{\modd}(A)}.
\]
Under the identification \eqref{AL-A}, $\widehat C^\bullet_{\underline{\modd}(A)}$ is characterized as the subring of $\widehat R^\bullet_{\underline{\modd}(A)}$
spanned by all homogeneous elements $g$ for which \eqref{diag-M-g} commutes with all irreducible $A$-modules $N$. Under that identification, 
$\widehat C^\bullet_{\underline{\modd}(A_{L})}$ is characterized as the subring of $\widehat R^\bullet_{\underline{\modd}(A_L)}$
spanned by all homogeneous elements $g$ for which \eqref{diag-M-g} commutes with the sub-collection (1) and (2) of $S$. 
(Here we apply \leref{comm-M-g} (b)). The commutativity condition on the elements of $\widehat C^\bullet_{\underline{\modd}(A_{L})}$
coming from the sub-collection (1) is precisely equivalent to the condition defining $\widehat C^\bullet_{\underline{\modd}(A)}$.
Therefore, a homogeneous element $g$ in $\widehat C^\bullet_{\underline{\modd}(A)}$ is in $\widehat C^\bullet_{\underline{\modd}(A_{L})}$ 
if and only if \eqref{diag-M-g} commutes for that $g$ and all $A_L$-modules of the form $M = \epsilon \otimes \kk_x$ for some $x \in L$. 
Each of the modules $\epsilon \otimes \kk_x$ is 1-dimensional and invertible: the inverse of $\epsilon \otimes \kk_x$ is $\epsilon \otimes \kk_{x^{-1}}$. 
These invertible objects form a group isomorphic to $L$, under $- \otimes -$.
By \leref{comm-M-g}(c), we may now conclude the following.

\bth{plavnik-witherspoon}
Let $A_{L}$ be the Plavnik-Witherspoon co-smash product Hopf algebras corresponding to the finite dimensional Hopf algebra $A$ and a finite group $L$, as defined above. Then there exist isomorphisms
\[ \widehat C^\bullet_{\underline{\modd}(A_{L})} \cong (\widehat C^\bullet_{\underline{\modd}(A)})^{L} 
\quad \mbox{and} \quad
C^\bullet_{\underline{\modd}(A_{L})} \cong (C^\bullet_{\underline{\modd}(A)})^L .\]
\eth

Many important families of finite dimensional Hopf algebras arise as Plavnik-Witherspoon co-smash product Hopf algebras, and thus satisfy 
\thref{plavnik-witherspoon}. 

\bex{twisted-quantum}[{\bf Twisted quantum Borels}] Retain the notation of \exref{QuantBorel}.
Denote by $\Gamma$ is the Dynkin quiver associated to ${\mathcal R}$. Its automorphism group $\Aut(\Gamma)$ acts on 
the lattices ${\mathbb Z}{\mathcal R}$ and $X$. If the lattice $\Lambda$ is stable under $\Aut(\Gamma)$, then $\Aut(\Gamma)$ 
acts on $u_{\zeta,\Lambda}({\mathfrak b})$ by Hopf algebra automorphisms by 
\[
\phi (E_{\alpha_i}) = E_{\phi(\alpha_i)}, \; \; 1 \leq i \leq n, \quad \phi (K_\mu) = K_{\phi(\mu)}, \; \; \mu \in \Lambda,
\]
for $\phi \in \Aut (\Gamma)$. 
One can consider the group algebra $k \Aut(\Gamma)$ and the co-smash product: 
\[
H:=(u_{\zeta,\Lambda}({\mathfrak b})^{*}\# k \Aut(\Gamma))^{*}.
\]

By \cite[Proposition 5.3.1(c)]{NVY2}, under the assumptions on $\zeta$ and ${\mathcal R}$, 
the categorical centers of the cohomology rings of the small quantum Borel subalgebras 
coincide with the full cohomology rings:
\[ 
\widehat C^\bullet_{\underline{\modd}(u_{\zeta,\Lambda}({\mathfrak b}))} = 
\widehat R^\bullet_{\underline{\modd}(u_{\zeta,\Lambda}({\mathfrak b}))} \cong \widehat \opH^\bullet(u_{\zeta,\Lambda}({\mathfrak b}), {\mathbb C}),
\; \; 
C^\bullet_{\underline{\modd}(u_{\zeta,\Lambda}({\mathfrak b}))}
= R^\bullet_{\underline{\modd}(u_{\zeta,\Lambda}({\mathfrak b}))} \cong \opH^\bullet(u_{\zeta,\Lambda}({\mathfrak b}), {\mathbb C}).
\]
Now, applying \exref{QuantBorel} we obtain
\[ 
\widehat C^\bullet_{\underline{\modd}(H)} \cong ( \widehat \opH^\bullet(u_{\zeta,\Lambda}({\mathfrak b}), {\mathbb C}) )^{\Aut(\Gamma)} 
\quad \mbox{and} \quad
C^\bullet_{\underline{\modd}(H)} \cong ( \opH^\bullet(u_{\zeta,\Lambda}({\mathfrak b}), {\mathbb C}) )^{\Aut(\Gamma)}.
\]
\eex

\bex{BW-ex}[{\bf Benson--Witherspoon Hopf algebras}] Let $G$ and $L$ be finite groups with $L$ acting on $G$ by group automorphisms, and let $k$ be a field of positive characteristic dividing the order of $G$. Consider the Benson-Witherspoon Hopf algebra $H_{G,L}$, which was studied in \cite{BW}. 
By definition, $H_{G,L}=(k[G] \# k L)^{*}$, where $k[G]$ is the dual of the group algebra of $G$, and $k L$ is the group algebra of $L$. 
\thref{plavnik-witherspoon} implies that
\[
C^\bullet_{\underline{\modd}(H_{G,L})} \cong \opH^\bullet(G, k)^L \quad \mbox{and} \quad
\widehat C^\bullet_{\underline{\modd}(H_{G,L})} \cong \widehat \opH^\bullet(G, k)^L. 
\]
\eex

\bex{Fin-Gr-Schemes-ex}[{\bf Finite group schemes}] Consider $k[\Omega]$, the coordinate algebra of $\Omega$, 
a finite group scheme, which is a commutative finite dimensional Hopf algebra. Assume for this example that $k$ is algebraically closed. 
One has $\Omega\cong \pi \ltimes \Omega_{0}$ where $\pi$ is a finite group and 
$\Omega_{0}$ is an infinitesimal group scheme \cite[Remark 4.1]{FP1}. Let $\mathrm{Dist}(\Omega)$ be the distribution algebra of $\Omega$. Then 
\begin{eqnarray*} 
k[\Omega]&\cong&  (\mathrm{Dist}(\Omega))^* 
                \cong (\mathrm{Dist}(\pi\ltimes \Omega_{0}))^{*} \\
                &\cong& (\mathrm{Dist}(\Omega_{0}) \# k\pi)^{*} 
                \cong (k[\Omega_{0}]^{*} \# k\pi)^{*}.
\end{eqnarray*} 
So, $k[\Omega]$
is a Plavnik-Witherspoon co-smash product Hopf algebra for $A=k[\Omega_{0}]$ and $L = \pi$.
One has as algebras $k [ \Omega] \cong k[\pi]\otimes k[\Omega_{0}]$ where $k[\pi]$ is isomorphic to the dual of the group algebra of $\pi$.
The coordinate algebra of $\Omega_{0}$, 
$k[\Omega_{0}]$, has an augmentation ideal, ${\mathcal I}$, that is nilpotent, i.e., $k[\Omega_{0}]$ is a unipotent Hopf algebra. 
\exref{UnipHopgalg} implies that
\[
\widehat C^\bullet_{\underline{\modd}(k[\Omega_0])} = \widehat R^\bullet_{\underline{\modd}( k[\Omega_{0}]) } 
\cong \widehat \opH^\bullet(k[\Omega_0], k)
\quad \mbox{and} \quad
C^\bullet_{\underline{\modd}(k[\Omega])} = R^\bullet_{\underline{\modd}( k[\Omega_{0}] )}
\cong \opH^\bullet(k[\Omega_0], k).
\]
\thref{plavnik-witherspoon} implies that
\[ 
\widehat C^\bullet_{\underline{\modd}(k[\Omega])} \cong ( \widehat \opH^\bullet(k[\Omega_0], k) )^{\pi} 
\quad \mbox{and} \quad
C^\bullet_{\underline{\modd}(k[\Omega])} \cong (\opH^\bullet(k[\Omega_0], k))^\pi .
\]
\eex
\section{From the Balmer spectrum to the central cohomological support}
In this section we prove part~(a) of Theorem B from the introduction, which produces continuous maps from the Balmer spectrum $\Spc \bK$ of an M$\Delta$C, $\bK$, to the 
target spaces of the Tate central cohomological support $\Spech \widehat C^\bullet$ and the target space of the central cohomological support $\Proj C^\bullet$.

%
\subsection{Properties of the central cohomological support maps and cones} We begin with the following proposition. 

\bpr{cent-supp}
For every M$\Delta$C, $\bK$, the Tate central cohomological support $\widehat  W_C$ satisfies the following properties:
\begin{enumerate}
\item[(a)] $\widehat W_C(0)= \varnothing$ and $\widehat W_C(\unit)= \Spech \widehat C^\bullet$;
\item[(b)] $\widehat W_C( A \oplus B) = \widehat W_C(A) \cup \widehat W_C(B)$ for all $A, B \in \bK$;
\item[(c)] $\widehat W_C( \Sigma A) = \widehat  W_C(A)$ for all $A \in \bK$;
\item[(d)] If $A _1\to A_2 \to A_3 \to \Sigma A_1$ is a distinguished triangle, then \\
$\widehat W_C(A_1) \subseteq \widehat  W_C(A_2) \cup \widehat  W_C(A_3)$;
\item[(e)] $\widehat W_C(A \otimes B) \subseteq \widehat W_C(A) \cap \widehat W_C(B)$ for all $A, B \in \bK$.
\end{enumerate}
The central support $W_C$ satisfies the same properties, with the second part of property (a) replaced by $W_C(\unit)= \Proj C^\bullet$.
\epr

This implies that the central cohomological support is a quasi-support datum, although (e) is a stronger condition than is required for quasi-support data.

\begin{proof} The proofs of properties (a)--(d) for both $\widehat  W_C$ and $W_C$ are analogous to the corresponding statements for the classical cohomological support
and left to the reader.

We prove (e) for the Tate central cohomological support. The argument for the central cohomological support is similar.
First we note that $\widehat W_C(A \otimes B) \subseteq \widehat W_C(A)$ is straightforward from the definition. Suppose $g: \unit \to \Sigma^i \unit$ in $\widehat C^\bullet$ is in $\Ann_{\widehat C^\bullet}(\hEnd^\bullet(A))$; in other words, $g \otimes \id_A= 0$, as maps $A \cong \unit \otimes A \to \Sigma^i \unit \otimes A \cong \Sigma^i A$. Then it is clear that $g \otimes \id_{A \otimes B}= g \otimes \id_A \otimes \id_B$ is also equal to 0, and so $g \in \Ann_{\widehat C^\bullet}(\hEnd^\bullet(A \otimes B)).$ Since 
$$\Ann_{\widehat C^\bullet}(\hEnd^\bullet(A)) \subseteq \Ann_{\widehat C^\bullet}(\hEnd^\bullet(A \otimes B)),$$
it follows that $\widehat W_C(A \otimes B) \subseteq \widehat W_C(A)$. 

For the inclusion $\widehat W_C(A \otimes B) \subseteq \widehat W_C(B)$, we must appeal to the definition of $\widehat C^\bullet$. The proof is an induction on the length of $A$ with respect to the generating set $\mc{K}$. In the case that $A$ is in $\mc{K}$, since $g \in \widehat C^\bullet$, the structure isomorphisms of $\bK$ identify $g \otimes \id_A \otimes \id_B$ with $\id_A \otimes g \otimes \id_B = \id_A \otimes 0 = 0$. 
Thus, $g \in \Ann_{\widehat C^\bullet} ( \hEnd^\bullet ( A \otimes B))$, and $\widehat W_C(A \otimes B) \subseteq \widehat W_C(B)$ follows as before. If $A$ is not in $\mc{K}$, we can pick a distinguished triangle
\[
A \to A_2 \to A_3 \to \Sigma A
\]
with the lengths of $A_2$ and $A_3$ less than the length of $A$ with respect to $\mc{K}$. Then 
\begin{align*}
\widehat W_C(A \otimes B) &\subseteq \widehat W_C(A_2 \otimes B ) \cup \widehat W_C(A_3 \otimes B)\\
&\subseteq (\widehat W_C(A_2) \cap \widehat W_C(B)) \cup (\widehat W_C(A_3) \cap \widehat W_C(B))\\
&\subseteq \widehat W_C(B).
\end{align*}
\end{proof}
\prref{cent-supp}(e) implies at once the following:
\bco{w-ideal}
For an M$\Delta$C, $\bK$, and $A \in \bK$,  
\[
\widehat W_C(A) = \Phi_{\widehat W_C}(\langle A \rangle) \quad \mbox{and} \quad W_C(A) = \Phi_{W_C}(\langle A \rangle). 
\]
\eco
Before constructing maps between the Balmer spectrum of a monoidal triangulated category and $\Spech \widehat C^\bullet$, $\Proj C^\bullet$, 
we recall an elementary consequence of the Octahedral Axiom for triangulated categories. Recall that given a morphism $g:A \to B$ in a triangulated category, $\cone(g)$ is the object (which is unique, up to a possibly non-unique isomorphism) given in a distinguished triangle
$$A \xrightarrow{g} B \rightarrow   \cone(g) \rightarrow  \Sigma A.$$

\ble{comp-cone-p}
Let $\bK$ be a triangulated category and $\bI$ any triangulated subcategory of $\bK$. Let $g': A \to B$ and $g: B \to C$ be morphisms in $\bK$. 
Suppose $\cone(g)$ and $\cone(g')$ are in $\bI$. Then $\cone(g g')$ is in $\bI$ as well. 
\ele
\begin{proof}
By the Octahedral Axiom for triangulated categories, there exist morphisms giving the following diagram, where the rows and columns are distinguished triangles:
\begin{center}
\begin{tikzcd}
A \arrow[r, "g'"] \arrow[d, Rightarrow, no head] & B \arrow[d, "g"] \arrow[r]                        & \cone(g') \arrow[r] \arrow[d]  & \Sigma A \arrow[d, Rightarrow, no head] \\
A \arrow[d] \arrow[r, "gg'"]                     & C \arrow[r] \arrow[d]                             & \cone(gg') \arrow[r] \arrow[d] & \Sigma A \arrow[d]                      \\
0 \arrow[d] \arrow[r]                            & \cone(g) \arrow[r, Rightarrow, no head] \arrow[d] & \cone(g) \arrow[r] \arrow[d]   & 0 \arrow[d]                             \\
\Sigma A \arrow[r]                               & \Sigma B \arrow[r]                                & \Sigma \cone(g') \arrow[r]     & \Sigma^2 A                             
\end{tikzcd}
\end{center} 
Since the third column is a triangle, if $\cone(g')$ and $\cone(g)$ are in $\bI$, so is $\cone(gg')$.
\end{proof}
\subsection{Construction of $\widehat \rho$ and $\rho$}
In the next two results we introduce a continuous map from the noncommutative Balmer spectrum to the homogeneous spectrum for the Tate categorical center. These results first appeared in the commutative case by Balmer in 
\cite[Definition 5.1, Theorem 5.3]{Balmer2}. 

We first recall an elementary property of $\cone(g)$, for $g \in \widehat R^\bullet$ (see \cite[Proposition 2.13]{Balmer2} in the symmetric case, or \cite[Proposition 3.6(d)]{BKSS}).

\bpr{f-tens-cone-f} 
Let $\bK$ be a M$\Delta$C, and 
\[
\unit \xrightarrow{g} \Sigma^i \unit \xrightarrow{h} \cone(g) \xrightarrow{f} \Sigma \unit
\]
be a distinguished triangle. Then 
\begin{enumerate}
\item[(a)] $g \otimes g \otimes \id_{\cone(g)}=0$; 
\item[(b)] $\id_{\cone(g)} \otimes g \otimes g =0$.
\end{enumerate}
\epr 

\begin{proof}
We first show (a). From the distinguished triangle given in the proposition, we obtain for any object $A$ a long exact sequence
\[
... \to \Hom^j (A, \unit) \xrightarrow{g \cdot} \Hom^j (A, \Sigma^i \unit) \to \Hom^j (A, \cone(g)) \to ...
\]
by basic properties of triangulated categories, see \cite[Lemma 1.1.10]{Neeman}, where the maps are obtained by composition with shifts of the morphisms $g,$ $h,$ and $f$. From this long exact sequence, we obtain the short exact sequence
\begin{equation}
\label{ses-gg}
0 \to \hHom^\bullet( A, \unit)/ ( g \cdot \hHom^\bullet (A, \unit))  \to \hHom^\bullet (A, \cone(g))
\end{equation}
\begin{equation*}
\to \ker(g \cdot |_{\hHom^\bullet (A, \unit)}) \to 0,
\end{equation*}
where we have used the fact that $\hHom^\bullet(A, \unit)[i] \cong \hHom^\bullet(A, \Sigma^i \unit)$ as graded vector spaces. In fact, this is a short exact sequence of $\widehat R^\bullet$-modules. From this short exact sequence, it is clear that the action of $g^2$ acting on $\hHom^\bullet(A, \cone(g))$ is 0 for any object $A$, and in particular this holds for $A= \cone(g)$. Hence, $g \otimes g \otimes \id_{\cone(g)}=0$. 

The proof of (b) follows similarly, but here we use an alternate action of $\widehat R^\bullet$ on $\hHom^\bullet(A,B)$, for arbitrary objects $A$ and $B$, given by setting $g .x$ equal to the composition
\[
A \cong A \otimes \unit \xrightarrow{x \otimes g} \Sigma^j B \otimes \Sigma^i \unit \cong \Sigma^{i+j} B,
\]
for $g: \unit \to \Sigma^i \unit$ and $x: A \to \Sigma^j B$. Note that the two actions coincide on $\hHom^\bullet(A,\unit)$, since the diagram
\begin{center}
\begin{tikzcd}
                                                              &                                                                                    & A \otimes \unit \arrow[r, "x \otimes g"]                                      & \Sigma^j \unit \otimes \Sigma^i \unit \arrow[rd, "\cong"]            &                                                              &                    \\
                                                              &                                                                                    & \Sigma^j \unit \otimes \unit \arrow[ru, "\id \otimes g"] \arrow[r, "\cong"]   & \Sigma^j(\unit \otimes \unit) \arrow[r, "\Sigma^j (\id \otimes g)"]  & \Sigma^j ( \unit \otimes \Sigma^i \unit) \arrow[rd, "\cong"] &                    \\
A \arrow[rruu, "\cong"] \arrow[rrdd, "\cong"'] \arrow[r, "x"] & \Sigma^j \unit \arrow[rrrr, "\Sigma^j g"] \arrow[ru, "\cong"] \arrow[rd, "\cong"'] &                                                                               &                                                                      &                                                              & \Sigma^{i+j} \unit \\
                                                              &                                                                                    & \unit \otimes \Sigma^j \unit \arrow[r, "\cong"'] \arrow[rd, "g \otimes \id"'] & \Sigma^j(\unit \otimes \unit) \arrow[r, "\Sigma^j (g \otimes \id)"'] & \Sigma^j(\Sigma^i \unit \otimes \unit) \arrow[ru, "\cong"']  &                    \\
                                                              &                                                                                    & \unit \otimes A \arrow[r, "g \otimes x"']                                     & \Sigma^i \unit \otimes \Sigma^j \unit \arrow[ru, "\cong"']           &                                                              &                   
\end{tikzcd}
\end{center}
commutes by the structure axioms of an M$\Delta$C. That is, $\widehat R^\bullet$ acts centrally on $\hHom^\bullet(A, \unit)$. Then the proof of (b) follows exactly the proof of (a), until the final step, when we note that the action of $g^2$ on $\id_{\cone(g)} \in \hHom^\bullet(\cone(g),\cone(g))$ is equal to $\id_{\cone(g)} \otimes g \otimes g.$
\end{proof}

\bpr{rho-hat-bullet-well-def}
For an M$\Delta$C, $\bK$, 
there is a well-defined map $\widehat \rho : \Spc \bK \to \Spech \widehat  C^\bullet$ given by
\[
\widehat \rho (\bP) =  \langle g \in \widehat C^\bullet: g \; \mbox{is homogeneous}, \; \cone (g) \not \in \bP \rangle
\]
for $\bP \in \Spc \bK$. 
\epr
\begin{proof} Fix $\bP \in \Spc \bK$. We will prove that $\widehat \rho(\bP)$ is a prime ideal of $\widehat  C^\bullet$.
By definition, $\widehat \rho(\bP)$ is an ideal; we will first show that the homogeneous elements of degree $i$ in $\widehat \rho(\bP)$ are precisely those 
elements of $\widehat C^\bullet$ of form $g: \unit \to \Sigma^i \unit$ such that $\cone(g) \not \in \bP$. 

We first show that this set is closed under addition. Suppose $g$ and $g'$ are two morphisms $\unit \to \Sigma^i \unit$ such that $\cone(g+g') \in \bP$. We must show that one of $\cone(g)$ and $\cone(g')$ is in $\bP$. 
By \prref{f-tens-cone-f}, we know that $g \otimes g \otimes \id_{\cone(g)}$ is the zero morphism 
\[
\unit \otimes \unit \otimes \cone(g) \cong \cone(g) \xrightarrow{0} \Sigma^i \unit \otimes \Sigma^i \unit \otimes \cone(g)\cong \Sigma^{2i}\cone(g),
\] 
and likewise for $g'$; similarly, $\id_{\cone(g)} \otimes g \otimes g$ and $\id_{\cone(g')} \otimes g' \otimes g'$ are also 0. Consider the morphism $\id_{\cone(g)} \otimes \id_A \otimes (g+g')^3 \otimes \id_{\cone(g')}$ 
for an object $A \in \mc{K}$, the generating set of $\bK$. This is a morphism 
\[
\cone(g) \otimes A \otimes \unit^{\otimes 3} \otimes \cone(g') \to \cone(g) \otimes A \otimes \Sigma^i \unit^{\otimes 3} \otimes \cone(g').
\] 
By exactness of the monoidal product, the cone of this morphism is 
\[
\cone(g) \otimes A \otimes \cone((g+g')^{ 3}) \otimes \cone(g').
\]
Since $\cone(g+g') \in \bP$ by assumption, so is $\cone((g+g')^{ 3})$, by \leref{comp-cone-p}. Now note that since $g$ and $g'$ are in $\widehat C^\bullet$ and $A$ is in $\mc{K}$, $g$ and $g'$ commute (up to isomorphism) with $\id_A$, and they skew-commute with each other by the graded-commutativity of $\widehat C^\bullet$. Since both $\id_{\cone(g)} \otimes g^2$ and $(g')^2 \otimes \id_{\cone(g')} =0$, this implies that the entire morphism $\id_{\cone(g)} \otimes \id_A \otimes (g+g')^3 \otimes \id_{\cone(g')} =0$. Hence, we have the following distinguished triangle, where (as noted above) the third term is in $\bP$:
\begin{align*}
\cone(g) \otimes A \otimes \cone(g') &\xrightarrow{0} \Sigma^{3i} \cone(g) \otimes A \otimes \cone(g')\\
& \to \cone(g) \otimes A \otimes \cone((g+g')^3) \otimes \cone(g') \to
\end{align*}
By the Splitting Lemma for triangulated categories (i.e., \cite[Corollary 1.2.5]{Neeman}), this implies that the third term of this triangle is isomorphic to 
\[
\big ( \Sigma^{3i} \cone(g) \otimes A \otimes \cone(g') \big ) \oplus \big ( \Sigma \cone(g) \otimes A \otimes \cone(g') \big ).
\]
By the thickness of $\bP$, each summand, and hence $\cone(g) \otimes A \otimes \cone(g')$, is in $\bP$. Since this holds for any object $A$ of $\mc{K}$, and $\mc{K}$ generates $\bK$ as a triangulated category, by induction $\cone(g) \otimes B \otimes \cone(g') \in \bP$ for any object $B$ of $\bK$, by the following argument. If $A_1$ and $A_2$ satisfy the property that $\cone(g) \otimes A_i \otimes \cone(g') \in \bP$ for $i=1,2$ and 
\[
A_1 \to A_2 \to A_3 \to \Sigma A_1
\] is a distinguished triangle, then 
\[
\cone(g) \otimes A_1 \otimes \cone(g') \to \cone(g) \otimes A_2 \otimes \cone(g') \to \cone(g) \otimes A_3 \otimes \cone(g')
\]
is also distinguished, and since $\bP$ is a triangulated subcategory, $\cone(g) \otimes A_3 \otimes \cone(g')$ is in $\bP$ as well. Similarly, if $A$ satisfies the property that $\cone(g) \otimes A \otimes \cone(g') \in \bP$, then 
\[
\cone(g) \otimes \Sigma A \otimes \cone(g') \cong \Sigma( \cone(g) \otimes A \otimes \cone(g'))
\] is in $\bP$ as well. Therefore $\cone(g) \otimes B \otimes \cone(g') \in \bP$ for all $B \in \bK$, and by the fact that $\bP$ is prime, either $\cone(g)$ or $\cone(g')$ is in $\bP$. This completes our first claim. 

Now suppose $g$ is one of the homogeneous generators $\widehat \rho(\bP)$, that is, $g: \unit \to \Sigma^i \unit$ for some $i$, and $\cone(g) \not \in \bP$. We now claim that given any other map $g' : \unit \to \Sigma^j \unit$ in $\widehat C^\bullet$, $gg'$ also satisfies $\cone(gg') \not \in \bP$. One way to see this is by Verdier localization: since $\bP$ is thick, there exists a triangulated category $\bK / \bP$ where objects are the same as in $\bK$, equipped with a triangulated functor $\bK \to \bK / \bP$ which is the identity on objects, and the objects which are isomorphic to 0 in $\bK / \bP$ are precisely those objects in $\bP$ (see \cite[Section 2.1]{Neeman}). If $\cone(gg')$ was in $\bP$, then the image of $gg'$ would be invertible in $\bK / \bP$, which would imply that both $g$ and $g'$ would have images also invertible in $\bK/\bP$, in other words, $\cone(g)$ and $\cone(g')$ would necessarily be in $\bP$. 

We have shown that the homogeneous generators of $\widehat \rho(\bP)$ are closed under addition and multiplication with arbitrary homogeneous elements of $\widehat C^\bullet$, and so the homogeneous elements of $\widehat \rho(\bP)$ are precisely those morphisms $g$ in $\widehat C^\bullet_i$ such that $\cone(g) \not \in \bP$. By \leref{comp-cone-p}, these homogeneous elements satisfy the prime condition, since if $\cone(gg') \not \in \bP$ then this implies that one of $\cone(g)$ and $\cone(g')$ must also be not in $\bP$. Thus, $\widehat \rho(\bP)$ is a homogeneous prime ideal of $\widehat C^\bullet$.
\end{proof}
\bpr{rho-hat-bullet-cont}
For all M$\Delta$Cs, $\bK$, the map $\widehat \rho: \Spc \bK \to \Spech \widehat  C^\bullet$ 
is continuous.
\epr
\begin{proof} For a proof in the symmetric case, see \cite[Theorem 5.3(b)]{Balmer2}. Consider an arbitrary closed set in $\Spech \widehat C^\bullet$; this has the form $Z(I) = \{ \mf{p} \in \Spech \widehat C^\bullet : I \subseteq \mf{p} \}$ for some homogeneous radical ideal $I$ of $\widehat C^\bullet$.
Define $\mc{S}$ as the collection of objects $\mc{S} = \{ \cone(g) : g \in I$ homogeneous$\}$. Now we note that for any $\bP$ in $\Spc \bK$, we have
\begin{align*}
\bP & \in (\widehat \rho)^{-1} ( Z(I))\\
& \Updownarrow\\
I & \subseteq \widehat \rho (\bP)\\
& \Updownarrow\\
\forall g \in I, g &\in \widehat \rho(\bP) \\
& \Updownarrow\\
\forall g \in I, \cone(g) &\not \in \bP\\
& \Updownarrow\\
\bP &\in  V(\mc{S}).
\end{align*}
Thus, the preimage of the closed set $Z(I)$ in $\Spech \widehat C^\bullet$ is the closed set $V (\mc{S})$ in $\Spc \bK$. 
\end{proof}

Denote the contraction map
\[
\theta : \Spech \widehat  C^\bullet \to \Spech  C^\bullet, \quad \mf{p} := \mf{p} \cap C^\bullet.
\]
Since this map is continuous, we have the following corollary of \prref{rho-hat-bullet-cont}:
\bco{rho-bullet-cont}
For an M$\Delta$C, $\bK$, the map 
\[
\rho:= \theta \circ \widehat \rho : \Spc \bK \to \Spech C^\bullet,
\]
which is explicitly given by
\[
\rho (\bP)  = \langle g \in \widehat C^\bullet: g \; \mbox{is homogeneous}, \; \cone (g) \not \in \bP \rangle \quad \mbox{for} \quad \bP \in \Spc \bK,
\]
is continuous.
\eco

Another consequence of \prref{rho-hat-bullet-cont} yields information about the inverse images of $\widehat{\rho}$ and $\rho$.

\bco{rho-comp} For an M$\Delta$C, $\bK$, 
\begin{enumerate}
\item[(a)] the map $\widehat \rho: \Spc \bK \to \Spech \widehat C^\bullet$ satisfies
\[
\widehat  \rho^{-1} (\widehat W_C(A)) \supseteq V(A) \quad \mbox{for all} \quad A \in \bK
\]
and
\item[(b)] the map $\rho: \Spc \bK \to \Spech C^\bullet$ satisfies
\[
\rho^{-1} (W_C(A)) \supseteq V(A) \quad \mbox{for all} \quad A \in \bK.
\]
\end{enumerate}
\eco

\begin{proof}
By definition, $\widehat W_C(A)=Z(I)$, where $I$ is the homogeneous ideal given as the annihilator of $\hEnd^\bullet(A)$ in $\widehat C^\bullet$. By the last part of the proof of \prref{rho-hat-bullet-cont}, $\widehat \rho^{-1}(\widehat W_C(A)) = V(\mc{S})$, where $\mc{S}$ is the collection of objects defined by
\[
\mc{S} = \{ \cone(g) : g \mbox{ homogeneous } g\in I\}.
\]
Suppose $g$ is a homogeneous element of $I$, i.e., $\cone(g) \in \mc{S}$. Then, by applying $- \otimes A$ to the distinguished triangle
\[
 \unit \xrightarrow{g} \Sigma^i \unit  \to \cone(g) \to \Sigma \unit,
\]
we obtain the distinguished triangle
\[
A \xrightarrow{0} \Sigma^i A \to  \cone(g) \otimes A \to \Sigma A.
\]
Hence, $\Sigma^i A\oplus \Sigma A \cong \cone(g) \otimes A$, and so $A$ is in the thick ideal generated by $\cone(g)$. In other words, if a prime ideal $\bP$ is in $V(A)$, meaning that it does not contain $A$, then it must not contain $\cone(g)$ either, and hence is in $V(\mc{S})$. This proves (a), and (b) follows directly. 
\end{proof}

\section{The surjectivity of $\rho$ and $\widehat \rho$}
\label{sec:surj}

We will show in this section that if $\bK$ is a M$\Delta$C satisfying the (wTfg) condition, then $\widehat \rho^\bullet$ is surjective, 
and that if $\bK$ satisfies the (wfg) condition then the image of $\rho^\bullet$ contains $\Proj C^\bullet_{\bK}$ (parts (b) and (d) of Theorem B in the introduction). 
In the important case of stable categories of finite tensor categories $\underline{\bT}$, we also prove part (c) of Theorem B that  
$\rho^\bullet$ takes values in $\Proj C^\bullet_{\underline{\bT}}$ (i.e., misses the irrelevant ideal of $\underline{\bT}$) 
and, with this codomain, $\rho^\bullet$ is surjective. 
\subsection{Behavior of cones under the central cohomological support maps}
We first recall that for $g \in \widehat C^\bullet$, the object $\cone(g)$ possesses desirable properties with respect to the central cohomological support. 
The objects $\cone(g)$ appearing here are analogues of Carlson's $L_{\zeta}$-modules, see \cite[Sect.\ 5.9]{Benson1}.

For homogeneous elements $g_1 \in \widehat C^\bullet$ (resp. $g_2 \in C^\bullet$), denote by
$\widehat Z(g_1)$ (resp. $Z(g_2)$) the Zariski closed subsets of $\Spech \widehat C^\bullet$  (resp. $\Spech C^\bullet$) defined by $g_1$ (resp. $g_2$).

The following result is given in \cite[Proposition 3.6, Proposition 3.7]{BKSS}. This is a monoidal triangulated version of \cite[Proposition 5.91]{Benson1}.

\bpr{lzeta-support} \cite{BKSS}
Let $\bK$ be an M$\Delta$C, $A$ an object of $\bK$, and $g: \unit \to \Sigma^i \unit$ a homogeneous element of $\widehat C^\bullet$. 
\begin{enumerate}
\item[(a)] $\widehat W_C(\cone(g) \otimes A) \subseteq \widehat W_C(A) \cap \{\mf{p} \in \Spech \widehat C^\bullet : g \in \mf{p}\}:= \widehat W_C(A) \cap \widehat Z (g).$
\item[(b)] If $\bK$ satisfies the (wTfg) condition, then $\widehat W_C(\cone(g) \otimes A) = \widehat W_C(A) \cap \widehat Z(g)$. 
\item[(c)] If $i \geq 0$, then $W_C(\cone(g) \otimes A) \subseteq W_C(A) \cap \{\mf{p} \in \Spech C^\bullet : g \in \mf{p}\}:= W_C(A) \cap Z(g).$
\item[(d)] If $i \geq 0$ and $\bK$ satisfies the (wfg) condition, then $W_C(\cone(g) \otimes A)$ and $W_C(A) \cap Z(g)$ coincide, except possibly at the unique maximal homogeneous ideal $\mf{m}$ of $C^\bullet$. 
\end{enumerate}
\epr

\begin{proof}
Since we have a distinguished triangle
\[
A \to  \Sigma^i A \to \cone(g) \otimes A \to \Sigma A,
\]
$\widehat W_C(\cone(g) \otimes A) \subseteq \widehat W_C(A)$. Moreover,  $\widehat W_C(\cone(g) \otimes A) \subseteq \widehat W_C(\cone(g))$, by the fact that the central support is a quasi-support datum. It remains to show that $\widehat W_C(\cone(g)) \subseteq \widehat Z(g)$. This follows from \prref{f-tens-cone-f}, since if $\mf{p}$ contains the annihilator of $\hEnd^\bullet(\cone(g))$ then it must contain $g^2$, and by primeness then contains $g$. This proves (a); (c) follows similarly.

For (b), suppose that $\mf{p}$ contains $g: \unit \to \Sigma^i \unit$, and the annihilator of $\hEnd^\bullet(A)$. We must show that $\mf{p} \in \widehat{W}_C(A \otimes \cone(g))$. Since the action of $\widehat C^\bullet$ on $\hHom^\bullet(B, A\otimes \cone(g))$ factors through $\hHom^\bullet(A \otimes \cone(g), A \otimes \cone(g))$ for any object $B$, it is enough to prove that 
$$\Ann_{\widehat C} \hHom^\bullet(B, A \otimes \cone(g)) \subseteq \mf{p}$$
for some $B$. Suppose to the contrary. From the distinguished triangle
\[
A \to A \to \cone(g) \otimes A \to \Sigma A,
\]
we obtain a short exact sequence of $\widehat C^\bullet$-modules by the same argument as for the short exact sequence (\ref{ses-gg}), for any object $B$:
\begin{equation}
\label{ses-gg-two}
0 \to \hHom^\bullet( B,A)/ ( g \cdot \hHom^\bullet (B,A))  \to \hHom^\bullet (B, \cone(g) \otimes A)
\end{equation}
\begin{equation*}
\to \ker(g \cdot |_{\hHom^\bullet (B, A)}) \to 0.
\end{equation*}
Since $\mf{p}$ does not contain the annihilator of $\hHom^\bullet(B,  \cone(g) \otimes A)$, we have that 
\[
\hHom^\bullet(B, \cone(g) \otimes A)_{\mf{p}}=0.
\] By (\ref{ses-gg-two}), this implies $\hHom^\bullet(B,A)_{\mf{p}} = g \cdot \hHom^\bullet(B,A)_{\mf{p}}$. Now using the hypothesis that $\hHom^\bullet( B,A)_{\mf{p}}$ is a finitely-generated $\widehat C_{\mf{p}}^\bullet$-module, and $g$ is in $\mf{p} \widehat C_{\mf{p}}$, by Nakayama's Lemma, $\hHom^\bullet(B,A)_{\mf{p}}=0$. Since $\hHom^\bullet(B,A)$ is a finitely-generated $\widehat C$-module by assumption, this implies that $\Ann_{\widehat C} \hHom^\bullet(B,A) \not \subseteq \mf{p}$, for any object $B$. But this is a contradiction, since we know that $\Ann_{\widehat C} \hEnd^\bullet(A) \subseteq \mf{p}$ by assumption. 

The proof of (d) follows similarly to (b), with one small modification. Suppose that $\mf{p} \in \Spech C^\bullet$ contains $g: \unit \to \Sigma^i \unit$, is different from the irrelevant (i.e. maximal homogeneous) ideal $\mf{m}$, and also contains the annihilator of $\End^\bullet(A).$ Suppose for the sake of contradiction that $\mf{p}$ does not contain the annihilator of $\Hom^\bullet(B, A \otimes \cone(g))$ for any object $B$. Using the long exact sequence as before, we obtain an analogue of (\ref{ses-gg-two}):
\begin{equation}
\label{ses-gg-three}
0 \to \Hom^{\geq i}( B,A)/ ( g \cdot \Hom^\bullet (B,A))  \to \Hom^\bullet (B, \cone(g) \otimes A)
\end{equation}
\begin{equation*}
\to \ker(g \cdot |_{\Hom^\bullet (B, A)}) \to 0.
\end{equation*}
This is a short exact sequence of $C^\bullet$-modules. Now since localization at $\mf{p}$ sends the middle term of this short exact sequence to 0 by assumption, we have 
\[
\Hom^{\geq i}(B,A)_{\mf{p}} = g \cdot \Hom^\bullet(B,A)_{\mf{p}}.
\] 
But now note that $\Hom^\bullet(B,A)_{\mf{p}} \cong \Hom^{\geq i} (B,A)_{\mf{p}}$ by commutative algebra, since $\mf{p}$ by assumption does not contain the irrelevant ideal of $C^\bullet$. This implies that there exists $t$ a positively-graded homogeneous element of $C^\bullet$ which is not contained in $\mf{p}$, and the map
\begin{align*}
\Hom^\bullet (B,A)_{\mf{p}} &\to \Hom^{\geq i}(B,A)_{\mf{p}},\\
\frac{x}{y} & \mapsto \frac{t^i x}{t^i y}
\end{align*}
is an isomorphism for $x \in \Hom^\bullet(B,A)$ and $y$ an element of $C^\bullet$ not in $\mf{p}$. Hence, $\Hom^\bullet(B,A)_{\mf{p}} = g \cdot \Hom^\bullet(B,A)_{\mf{p}}$, and the remainder of the proof of (d) follows exactly the end of the proof of (b). 
\end{proof}

The next three corollaries follow immediately from \prref{lzeta-support}.

\bco{m-in-im}
Suppose $\bK$ is an M$\Delta$C, and that there exists a finitely generated ideal $I$ such that $\sqrt{I}$ is the unique maximal homogeneous ideal $\mf{m}$ of $C^\bullet$. If $\mf{m}$ is in the image of $\rho: \Spc \bK \to \Spech C^\bullet$, then there exists an object $A$ of $\bK$ such that $W_C(A) = \{ \mf{m}\}$.
\eco

\begin{proof}
Suppose $\rho(\bP) = \mf{m}$, and let $I= \langle g_1,..., g_n\rangle$. Since $\rho(\bP)=\mf{m}$, $\cone(g_i) \not \in \bP$ for all $i$. This implies by primeness of $\bP$ and the fact that each $\cone(g_i)$ $\otimes$-commutes with all objects of $\bK$ that $\cone(g_1) \otimes... \otimes \cone(g_n) \not = 0$. By \prref{lzeta-support}(c), 
\[
W_C(\cone(g_1) \otimes... \otimes \cone(g_n)) \subseteq Z(g_1,..., g_n) = \{ \mf{m}\}. 
\]
Since $\cone(g_1) \otimes... \otimes \cone(g_n)$ is nonzero, $W_C(\cone(g_1) \otimes... \otimes \cone(g_n)) \not = \varnothing$, and so 
\[
W_C(\cone(g_1) \otimes... \otimes \cone(g_n) ) = \{ \mf{m}\}. 
\]
\end{proof} 

\bco{stmod-fintens-max}
Let $\bT$ be a finite tensor category and $\bK = \underline{\bT}$ be its stable category. Suppose that $\bK$ satisfies (wfg). Then $\mf{m}$, 
the irrelevant ideal of $C^\bullet$, is not in the image of $\rho: \Spc \bK \to \Spech C^\bullet$, i.e., $\rho$ takes values in $\Proj C^\bullet$:
\[
\rho: \Spc \bK \to \Proj C^\bullet.
\]
\eco

\begin{proof}
This follows from \cite[Corollary 4.2]{BPW1}, which states that the dimension of the support variety of an object $A$ is equal to 0 if and only if $A$ is projective in $\bT$ (see also \cite[Proposition 2.4(1)]{FW1} in the finite-dimensional Hopf algebra setting). While the support varieties of \cite{BPW1} are based on the whole cohomology ring $R^\bullet$ rather than its categorical center $C^\bullet$, the same proofs carry through assuming the finite generation of each $\End^\bullet(A)$ over $C^\bullet$, rather than its finite generation over $R^\bullet$. Hence, there is no object $A$ such that $W_C(A) =\{\mf{m}\}$, since this would imply that $A$ is projective in $\bT$, and thus $W_C(A) = W_C(0) = \varnothing$. By \coref{m-in-im}, $\mf{m}$ is not in the image of $\rho$. 
\end{proof}

\bco{tensor-zero}
Let $\bK$ be an M$\Delta$C.
\begin{enumerate}
\item[(a)] Suppose $\bK$ satisfies (wTfg). Let $g_1,..., g_n \in \widehat C^\bullet$ and $\cone(g_1) \otimes ... \otimes \cone(g_n) \cong 0$. Then $\langle g_1,..., g_n \rangle = \widehat C^\bullet$. 
\item[(b)] Suppose $\bK$ satisfies (wfg). Let $g_1,..., g_n \in C^\bullet$ and $\cone(g_1) \otimes... \otimes \cone(g_n) \cong 0$. Then $\mf{m} \subseteq \sqrt{\langle g_1,..., g_n \rangle}$, where $\mf{m}$ denotes the irrelevant ideal of $C^\bullet$. If $g_i$ is of strictly positive degree for $i=1,2, \dots, n$, then $ \mf{m} = \sqrt{\langle g_1,..., g_n \rangle}.$
\end{enumerate}
\eco

\begin{proof}
By \prref{lzeta-support}(b), we know that
\begin{align*}
\widehat W_C( \cone(g_1) \otimes... \otimes \cone(g_n)) &= \bigcap_i \widehat W_C(\cone(g_i))\\
&= \{ \mf{p} \in \Spech \widehat C^\bullet: g_i \in \mf{p} \; \forall \; i\}.
\end{align*}
Clearly, if $\cone(g_1) \otimes... \otimes \cone(g_n) \cong 0$, then from the definition of support, $\widehat W_C( \cone(g_1) \otimes... \otimes \cone(g_n)) = \varnothing$, which implies that $$\{ \mf{p} \in \Spech \widehat C^\bullet: g_i \in \mf{p}\; \forall \; i\} = \varnothing,$$ in other words, the ideal generated by $g_1,..., g_n$ is the entire ring $\widehat C^\bullet$. This proves (a).

Similarly to the proof of part (a), if $\cone(g_1) \otimes ... \otimes \cone(g_n) \cong 0$ for $g_i$ homogeneous in $C^\bullet$, then $W_C( \cone(g_1) \otimes... \otimes \cone(g_n)) = \varnothing$. This coincides with 
\[
Z(g_1,..., g_n) = \{ \mf{p} \in \Spech C^\bullet: g_i \in \mf{p} \; \forall \; i\} 
\]
except possibly for the maximal ideal $\mf{m}$, by \prref{lzeta-support}(d). Hence, $Z(g_1,..., g_n)$ is either empty, in which case $g_1,..., g_n$ generate the entire ring $C^\bullet$, or it is $\{\mf{m}\}$, in which case the radical of $\langle g_1,..., g_n \rangle$ is equal to $\mf{m}$. This proves (b). 
\end{proof}
\subsection{Surjectivity of  $\rho$ and $\widehat \rho$}

We can now show that finite generation insures that  $\rho$ and $\widehat \rho$ are surjective. Note that there is no requirement that $\widehat W_C$ or $W_C$ be weak support data, which is the requirement which ensures the existence of the maps $\eta_{\widehat W}: \Spech \widehat C^\bullet \to \Spc \bK$ and $\eta_{W}: \Proj C^\bullet \to \Spc \bK$, as in \thref{weak-univ}. 

\bth{rho-surj}
Let $\bK$ be an M$\Delta$C.
\begin{enumerate}
\item[(a)] If $\bK$ satisfies (wTfg), then $\widehat \rho: \Spc \bK \to \Spech \widehat C^\bullet$ is surjective.
\item[(b)] If $\bK$ satisfies (wfg) and $\mf{p}$ is a homogeneous prime of $C^\bullet$ different from the irrelevant ideal $\mf{m}$, then $\mf{p}$ is in the image of $\rho: \Spc \bK \to \Spech C^\bullet$.
\end{enumerate}
\eth

\begin{proof}
We prove (a) directly, and (b) follows from the analogous argument. Let $\mf{p}$ be a homogeneous prime in $\widehat C^\bullet$. We will construct a prime ideal $\bP$ in $\Spc \bK$ such that $\widehat \rho^\bullet(\bP)= \mf{p}$. To construct this prime ideal, we will employ \thref{maximal}.

Consider the following subsets of $\bK$:
\begin{enumerate}
\item $\bI= \{ A \in \bK : \mf{p} \not \in \widehat W_C(A) \}$.
\item $\mc{M} =\{ \cone(g_1) \otimes \cone(g_2) \otimes... \otimes \cone(g_n) : g_i \in \mf{p}$ homogeneous$\}.$
\end{enumerate}

First, note that $\bI$ is a thick ideal. This follows directly from the fact that $\widehat W_C$ is a quasi-support datum with the additional property that $\widehat W_C(A \otimes B) \subseteq \widehat W_C(A) \cap \widehat W_C(B)$. 

Second, note that $\mc{M}$ is closed under the tensor product, by definition. We claim that $\mc{M}$ is a multiplicative subset; it remains to be shown that all objects of $\mc{M}$ are not isomorphic to 0. This follows from \coref{tensor-zero}(a), since the ideal generated by $g_1,..., g_n$ is contained in $\mf{p}$ for any collection $g_1,..., g_n \in \mf{p}$. 

Next, we claim that $\bI$ and $\mc{M}$ are disjoint. Take an arbitrary element $\cone(g_1) \otimes... \otimes \cone(g_n)$ of $\mc{M}$. By \prref{lzeta-support}(b), 
\[
\widehat W_C(\cone(g_1) \otimes... \otimes \cone(g_n)) = \bigcap_i \widehat W_C(\cone(g_i)).
\]
Since each $g_i \in \mf{p}$, we have for all $i$ that $\mf{p} \in \widehat W_C(\cone(g_i))$, again by \prref{lzeta-support}. Hence, $\mf{p} \in \bigcap_i \widehat  W_C(\cone(g_i)) = \widehat W_C(\cone(g_1) \otimes ... \otimes \cone(g_n))$, which implies that $\cone(g_1) \otimes... \otimes \cone(g_n) \not \in \bI$. Thus, $\bI$ and $\mc{M}$ are disjoint.

Now, by \thref{maximal} we know that there exists at least one prime ideal $\bP$ of $\bK$ such that $\bI \subseteq \bP$ and $\mc{M} \cap \bP = \varnothing$. We claim that $\widehat \rho( \bP ) = \mf{p}$. By definition,
\[
\widehat \rho( \bP) = \{ g \in \widehat C^\bullet : \cone(g) \not \in \bP \}.
\] 
If $g \in \mf{p}$, then since $\bP \cap \mc{M} = \varnothing$, $\cone(g) \not \in \bP$, and hence $g \in \widehat \rho(\bP)$. On the other hand, if $h \not \in \mf{p}$, then $\mf{p} \not \in \widehat W_C(\cone(h))$ by \prref{lzeta-support}(b), and so $\cone(h) \in \bI$.  Since $\bI \subseteq \bP$, it follows that $\cone(h) \in \bP$, and hence $h \not \in \rho (\bP)$. Therefore, $\rho^\bullet(\bP) = \mf{p}$. 
\end{proof}

\section{Inverses of the maps $\rho$ and $\widehat \rho$}
\label{sec:final-thm}
In this section we prove that under natural conditions (which are stronger than those used in the previous sections),
the pairs $(\rho, \eta)$ and $(\widehat \rho, \widehat \eta)$ are pairs of inverse homeomorphisms. 
Under those conditions we obtain a classification of the thick two-sided ideals of an M$\Delta$C
in terms of homological data.
\subsection{Right inverses of $\rho$ and $\widehat \rho$}
If the central cohomological support $W_C : \bK \to \Proj C^\bullet$ is a weak support datum, then by \thref{weak-univ}
there is a unique continuous map $\eta: \Proj C^\bullet \to \Spc \bK$ defined by 
$$
\eta (\mf{p}):= \{ M \in \bK : \mf{p} \not \in \Phi_{W_C} ( \langle M \rangle) \} \quad \mbox{for} \quad \mf{p} \in \Proj C^\bullet.
$$
By \coref{w-ideal}, $$\eta(\mf{p}) = \{ M \in \bK : \mf{p} \not \in W_C(M) \}.$$

If the Tate central cohomological support $\widehat W_C : \bK \to \Spech \widehat C^\bullet$ is a weak support datum, then again 
by \thref{weak-univ} there is a unique continuous map $\widehat \eta: \Spech \widehat C^\bullet \to \Spc \bK$ given by 
$$
\widehat \eta(\mf{p}) := \{ M \in \bK : \mf{p} \not \in \Phi_{\widehat W_C} ( \langle M \rangle) \} \quad \mbox{for} \quad \mf{p} \in \Spech \widehat C^\bullet
$$
and by \coref{w-ideal}, $\widehat \eta(\mf{p}) = \{ M \in \bK : \mf{p} \not \in \widehat W_C (M) \}.$
\bpr{rho-f-contained}
Let $\bK$ be an M$\Delta$C for which the central cohomological support is a weak support datum. Then the following hold: 
\begin{enumerate}
\item[(a)] $\rho( \eta( \mf{p})) \subseteq \mf{p}$ for every homogeneous prime ideal $\mf{p}$ of $C^\bullet$.
\item[(b)] If $\bK$ satisfies the (wfg) condition, then $\rho( \eta( \mf{p})) = \mf{p}$ for every $\mf{p} \in \Proj C^\bullet$.
\item[(c)] If $\bK$ satisfies the (wfg) condition, then $\eta( \rho ( \bP)) \subseteq \bP$ for every prime ideal $\bP$ of $\bK$.
\end{enumerate}
\epr
\begin{proof}
(a) Let $\mf{p}$ be a homogeneous prime of $C^\bullet$. Then
\begin{align*}
\rho (\eta( \mf{p}))&= \langle g \in C^\bullet : \cone(g) \not \in \{M  \in \bK : \mf{p } \not \in \Phi_{W_C}( \langle M \rangle )\} \rangle\\
&= \langle g \in C^\bullet : \mf{p} \in \Phi_{W_C}(  \langle \cone (g) \rangle) \rangle\\
&= \langle g \in C^\bullet : \mf{p} \in W_C(\cone (g)) \rangle\\
&= \langle g \in C^\bullet: \Ann_{C^\bullet} ( \End^\bullet(  \cone(g))) \subseteq \mf{p} \rangle.
\end{align*}
Suppose that $g$ is one of the homogeneous generators of $\rho( \eta( \mf{p}))$. Then note that $g \otimes g$ is in the annihilator of $\cone(g)$, since $g \otimes g \otimes \id_{\cone(g)} = 0$, by \prref{f-tens-cone-f}. Hence, $g \otimes g \in \mf{p}$, and thus $g \in \mf{p}$. Since $\mf{p}$ contains all homogeneous generators of 
$\rho^\bullet ( \eta( \mf{p}))$, it contains the entire ideal.

(b) One has 
$$\rho( \eta( \mf{p})) = \langle g \in C^\bullet: \mf{p} \in W_C(\cone(g) ) \rangle.$$ 
By \prref{lzeta-support}(d), applied for $A = \unit$, $W_C(\cone(g)) = Z(g)$ or  $W_C(\cone(g)) \sqcup \{ \mf{m} \} = Z(g),$ where $\mf{m}$ is the irrelevant ideal of $C^\bullet$. 
Since $\mf{p} \neq \mf{m}$,
$$\rho( \eta( \mf{p})) = \langle g \in C^\bullet: \mf{p} \in Z(g) \rangle = \langle g \in C^\bullet: g \in \mf{p} \rangle = \mf{p}.$$ 

(c) One has 
\begin{align*}
\eta( \rho(\bP)) &= \{ M : \langle g: \cone(g) \not \in \bP  \rangle \not \in W_C(M) \}\\
&= \{ M : \Ann_{C^{\bullet}}( \End^\bullet(M)) \not \subseteq \langle g: \cone(g)  \not \in \bP \rangle \}.
\end{align*}
Suppose $M$ is in $\eta( \rho(\bP))$. Then there exists some $h \in \Ann_{C^{\bullet}}(\End^\bullet(M))$, i.e., having $h \otimes \id_M =0$, such that $h \not \in \langle g: \cone(g) \not \in \bP\rangle$; in particular, $h$ is not one of the generators of this ideal, and so $\cone(h) \in \bP$. Since $h \otimes \id_M =0$, the triangle
$$M \xrightarrow{0} \Sigma^i M  \to M \otimes \cone(h) \to \Sigma M$$ is distinguished. By the Splitting Lemma for triangulated categories, this implies that
$$M \otimes \cone(h) \cong \Sigma^i M \oplus \Sigma M.$$ The left hand side is in $\bP$ by the ideal property, and this implies by the thickness of $\bP$ that $M \in \bP$. This proves (c). 
\end{proof}
\prref{rho-f-contained}(b) and \coref{stmod-fintens-max} imply the following:
\bco{right-inv-rho} Assume that $\bT$ is a finite tensor category whose stable category $\underline{\bT}$ satisfies (wfg) and that the central cohomological support is a weak support datum.  
Then the map $\eta : \Proj C^\bullet \to \Spc \bK$ is a right inverse of $\rho : \Spc \bK \to \Proj C^\bullet$.
\eco
In a similar way to the proof of \prref{rho-f-contained}, by using \prref{lzeta-support}(b), one proves the following:
\bpr{rho-f-contained2}
Let $\bK$ be an M$\Delta$C for which the Tate central cohomological support is a weak support datum. Then the following hold: 
\begin{enumerate}
\item[(a)] $\widehat \rho( \widehat \eta( \mf{p})) \subseteq \mf{p}$ for every homogeneous prime ideal $\mf{p}$ of $\widehat C^\bullet$.
\item[(b)] If $\bK$ satisfies the (wTfg) condition, then $\widehat \rho( \widehat \eta( \mf{p})) = \mf{p}$ for every homogeneous prime ideal $\mf{p}$ of $\widehat C^\bullet$, 
i.e., $\widehat \eta$ is a right inverse of $\widehat \rho$.
\item[(c)] If $\bK$ satisfies the (wTfg) condition, then $\widehat \eta( \widehat \rho ( \bP)) \subseteq \bP$ for every prime ideal $\bP$ of $\bK$.
\end{enumerate}
\epr
\subsection{Conditions for $\rho$ and $\eta$ to be inverse homeomorphisms}  The following theorems in this section provide the description of the Balmer spectrum and the classification of thick tensor ideals for the compact part of an arbitrary compactly generated M$\Delta$C via the central cohomology rings. 
\bth{rho-inv-f} Let $\bK$ be an M$\Delta$C which is the compact part of a compactly generated  M$\Delta$C, $\widetilde{\bK}$.
Assume that $\bK$ satisfies the (wfg) condition, $\Proj C^\bullet$ is a Zariski space, and 
the central cohomological support of $\bK$ has an extension to a faithful extended weak support datum $\widetilde{\bK} \to \mc{X}(\Proj C^\bullet)$.
Then the following hold: 
\begin{enumerate}
\item[(a)] The maps $\eta$ and $\rho$ are inverse homeomorphisms
\[
\Spc \bK \mathrel{\mathop{\rightleftarrows}^{\rho}_\eta} \Proj C^\bullet.
\]
\item[(b)] The map 
\[
\Phi_{W_C} : \ThickId(\bK) \to \mc{X}_{sp}(\Proj C^\bullet)
\]
is an isomorphism of ordered monoids, where the set of thick ideals of $\bK$ is equipped
with the operation $\bI, \bJ \mapsto \langle \bI \otimes \bJ \rangle$ and the inclusion partial order, 
and $\mc{X}_{sp}(\Proj C^\bullet)$ is equipped with the operation of intersection and the inclusion partial order.
\end{enumerate}
\eth 
\begin{proof} Since $\Proj C^\bullet$ is a Zariski space, every closed set of $\Proj C^\bullet$ is the variety generated by a finite collection of homogeneous elements, 
say $g_1,..., g_n$. By \prref{lzeta-support}(d),  
\begin{align*}
W_C( \cone(g_1) \otimes... \otimes \cone(g_n)) &= W_C (\cone(g_1)) \cap... \cap W_C(\cone(g_n))\\
& = Z(g_1) \cap ... \cap Z(g_n)\\
&= Z(g_1,..., g_n),
\end{align*}
where for an element $g \in C^\bullet$, $Z(g)$ denotes the closed set of $\Proj C^\bullet$ defined by $g$ (and not of $\Spech C^\bullet$ as used in 
Section \ref{sec:surj}).
Hence, the central support satisfies the realization condition. Therefore, the extension of the central cohomological support
satisfies the assumptions of \thref{reconstr}. Part (a) of that theorem gives that $\eta$ is a homeomorphism. 
In view of this, \prref{rho-f-contained}(b) implies that $\rho$ takes values in $\Proj C^\bullet$ and is a right inverse of $\eta$. 
Hence, $\rho$ is an inverse of $\eta$, which proves part (a). Part(b) follows from \thref{reconstr}(b).
\end{proof}
In a similar manner, by using Propositions \ref{plzeta-support}(b) and \ref{prho-f-contained2}(b) and \thref{reconstr}, one proves the following:
\bth{rho-inv-hat-f} Let $\bK$ be an M$\Delta$C, which is the compact part of a compactly generated  M$\Delta$C, $\widetilde{\bK}$.
Assume that $\bK$ satisfies the (wTfg) condition, $\Spech \widehat C^\bullet$ is a Zariski space and 
the Tate central cohomological support of $\bK$ has an extension to a faithful extended weak support datum $\widetilde{\bK} \to \mc{X}(\Spech \widehat C^\bullet)$.
Then the following hold: 
\begin{enumerate}
\item[(a)] The maps $\widehat \eta$ and $\widehat \rho$ are inverse homeomorphisms
\[
\Spc \bK \mathrel{\mathop{\rightleftarrows}^{\widehat \rho}_{\widehat \eta}} \Spech \widehat C^\bullet.
\]
\item[(b)] The map 
\[
\Phi_{\widehat W_C} : \ThickId(\bK) \to \mc{X}_{sp}(\Spech \widehat C^\bullet)
\]
is an isomorphism of ordered monoids, where the set of thick ideals of $\bK$ is equipped
with the operation $\bI, \bJ \mapsto \langle \bI \otimes \bJ \rangle$ and the inclusion partial order, 
and $\mc{X}_{sp}(\Spech \widehat C^\bullet)$ is equipped with the operation of intersection and the inclusion partial order.
\end{enumerate}
\eth 

\section{Examples to Theorem D} 
\label{sec:ex-ThmD}
In this section we show that the assumptions of Theorem D can be 
verified in a uniform fashion for wide classes of Hopf algebras, leading to a classification of 
the thick tensor ideals of their stable module categories in terms of homological data
that goes far beyond previous classes of Hopf algebras treated on a case by case basis.

\subsection{Plavnik-Witherspoon co-smash product Hopf algebras} In the case of
Plavnik-Witherspoon co-smash products $A_L$ considered in Section \ref{plav-with}, 
one can prove that, if the initial Hopf algebra $A$ satisfies the conditions in \thref{rho-inv-f}, 
then so does $A_L$. This leads to a classification of the thick tensor ideals and noncommutative Balmer spectrum 
of the stable module category of $A_{L}$.

\bth{plavnik-witherspoon-classify} Let $A_{L}$ be the Plavnik-Witherspoon co-smash product Hopf algebra 
corresponding to the finite dimensional Hopf algebra $A$, and the finite group $L$ acting on $A$. 
Assume that the Hopf algebra $A$  satisfies the following conditions:
\begin{itemize} 
\item[(a)] $\underline{\modd}(A)$ satisfies the (wfg) condition;
\item[(b)] $C^{\bullet}_{\underline{\modd}(A)}$ is finitely generated algebra over the base field $k$;
\item[(c)] the central cohomological support of $\underline{\modd}(A)$ has an extension to a faithful extended weak support datum 
$\underline{\Mod}(A) \to \mc{X}(\Proj C_{\underline{\modd}(A)}^\bullet)$.
\end{itemize} 
Then $\underline{\modd}(A_L)$ also satisfies conditions (a)-(c). In particular,
\[
\Spc \big( \underline{\modd}(A_{L}) \big) \mathrel{\mathop{\rightleftarrows}^{\rho}_\eta} \Proj (C_{\underline{\modd}(A)}^\bullet)^L
\] 
are inverse homeomorphisms and the map 
\[
\Phi_{W_C} : \ThickId \big(\underline{\modd}(A_{L}) \big) \to \mc{X}_{sp}\big(\Proj (C_{\underline{\modd}(A)}^\bullet)^L \big)
\]
is an isomorphism of ordered monoids. 
\eth

\begin{proof} The theorem will follow from \thref{rho-inv-f} by showing that conditions (a), (b) and (c) hold when $A$ is replaced by $A_{L}$. 
By \thref{plavnik-witherspoon},
\[
C^\bullet_{\underline{\modd}(A_{L})} \cong (C^\bullet_{\underline{\modd}(A)})^L,
\]
and thus condition (b) holds for $A_L$ since $L$ is a finite group. Consequently, one has the geometric quotient map
\[
\nu : \Proj C^\bullet_{\underline{\modd}(A)} \to ( \Proj C^\bullet_{\underline{\modd}(A)} )/L \cong \Proj (C^\bullet_{\underline{\modd}(A)})^L
\cong \Proj (C^\bullet_{\underline{\modd}(A_L)}).
\]
To verify condition (a) for $A_L$ we apply the decomposition \eqref{AL-decomp}. It suffices to show that 
$\Ext^{\bullet}_{A_{L}}(M_x \otimes \kk_x,M_y\otimes \kk_y)$ is finitely generated  
$C^\bullet_{\underline{\modd}(A_L)}$-module for all $x,y \in L$, $M_x, M_y \in {\underline{\modd}}(A)$. 
This follows by observing that 
\begin{equation}
\label{Ext-isom}
\Ext^{\bullet}_{ {\underline{\modd}}(A_{L})}(M_x \otimes \kk_x,M_y\otimes \kk_y)\cong \Ext^{\bullet}_{ {\underline{\modd}}(A)}(M_x,M_y)\otimes \Hom_{k[L]}(\kk_x, \kk_y),
\end{equation}
which is finitely generated as a $C^\bullet_{\underline{\modd}(A_L)}$-module by the condition (a) for $A$ and the fact that 
$C^\bullet_{\underline{\modd}(A)}$ is module finite over $(C^\bullet_{\underline{\modd}(A)})^L \cong C^\bullet_{\underline{\modd}(A_L)}$.

Denote by $W_{C,A}$ and $W_{C,A_L}$ the central cohomological supports of ${\underline{\modd}(A)}$ and ${\underline{\modd}(A_L)}$, respectively.
Denote by $\widetilde{W}_{C,A}$ the extension of $W_{C,A}$ from condition (c) to ${\underline{\Mod}(A)}$. 
Let $\mathrm{For} : \underline{\Mod}(A_L) \to \underline{\Mod}(A)$ be
the forgetful functor associated to the canonical embedding of $A$ in $A_L$. It restricts to a functor $\mathrm{For} : \underline{\modd}(A_L) \to \underline{\modd}(A)$.
It follows from \eqref{Ext-isom} that
\[
W_{C,A_L}(M) = \nu (W_{C,A}( \mathrm{For}(M))) \quad \mbox{for all} \quad M \in {\underline{\modd}(A)}.
\]
Define the map 
\[
\widetilde{W}_{C,A_L} : \underline{\Mod}(A) \to \mc{X}(\Proj C_{\underline{\modd}(A)}^\bullet) 
\quad \mbox{by} \quad 
\widetilde{W}_{C,A_L} (M) := \nu (W_{C,A_L}( \mathrm{For}(M)))
\]
for $M \in {\underline{\Mod}(A)}$. Using the tensor product formula \cite[Theorem 2.3]{PW} for $A_L$-modules, one easily verifies that 
$\widetilde{W}_{C,A_L}$ is an extended weak support datum. Its faithfulness is obvious.
\end{proof}

\subsection{Benson-Witherspoon Hopf algebras} Consider a Benson-Witherspoon Hopf algebra $H_{G,L}=(k[G] \# k L)^{*}$
as in \exref{BW-ex}. It is a special case of the construction in the previous section with $A= k G$.
Conditions (a) and (b) of 
\thref{plavnik-witherspoon-classify} hold because $A$ is the group algebra of a finite group. 
A faithful extension of the cohomological support of $k G$ was constructed in \cite{BCR1,BCR2}, which verifies 
condition (c) in the theorem. By \thref{plavnik-witherspoon-classify} and \exref{BW-ex}, we have: 
\bpr{BW-class} For each Benson-Witherspoon Hopf algebra $H_{G,L}$, 
\[
\Spc \underline{\modd}(H_{G,L}) \mathrel{\mathop{\rightleftarrows}^{\rho}_\eta} \Proj (\opH^\bullet(G,k))^{L}
\]
are inverse homeomorphisms and the map 
\[
\Phi_{W_C} : \ThickId \big(\underline{\modd}(H_{G,L}) \big) \to \mc{X}_{sp}\big(\Proj (\opH^\bullet(G,k))^{L} \big)
\]
is an isomorphism of ordered monoids. 
\epr
A classification of the thick tensor ideals and the noncommutative Balmer spectra of 
the stable module categories of the Benson-Witherspoon Hopf algebras was given in 
\cite[Theorem 9.3.2]{NVY1} in terms homogeneous prime ideals of $(\opH^\bullet(G,k))^{L}$. 

\subsection{Finite group schemes} Let $k[\Omega]$ be the coordinate algebra of a finite group scheme $\Omega$ with $k$ algebraically closed.  Recall that 
$$k[\Omega]\cong  (k[\Omega_{0}] ^{*}\# k\pi)^{*},$$ 
where $\pi$ is a finite group and $\Omega_{0}$ is an infinitesimal group scheme. 
This is a Plavnik-Witherspoon co-smash product Hopf algebra for $A=k[\Omega_{0}]$ and $L = \pi$.
Since we have the isomorphisms of algebras
$$
A=k[\Omega_{0}] \cong k[t_{1},t_{2},\dots, t_{s}]/(t_{1}^{p^{r_{1}}},\dots, t_{s}^{p^{r_{s}}})
\cong
k({\mathbb Z}_{p^{r_{1}}}\times \dots \times {\mathbb Z}_{p^{r_{s}}})
$$ 
and $C_{\underline{\modd}(k[\Omega_{0}])}^\bullet \cong H^{\bullet}(k[\Omega_{0}],k)$, 
conditions (a) and (b) of \thref{plavnik-witherspoon-classify} hold. 
Here we use that, since $k[\Omega_{0}]$ is a unipotent Hopf algebra, its cohomological support 
depends only on its algebra structure, but not coalgebra structure.  

Next we verify condition (c) of \thref{plavnik-witherspoon-classify} for $k[\Omega_{0}]$. 
Let $E$ be the unique elementary abelian subgroup of
${\mathbb Z}_{p^{r_{1}}}\times \dots \times {\mathbb Z}_{p^{r_{s}}}$ consisting of the identity element and elements of order $p$. 
One can now apply the constructions in \cite[Section 2]{BCR2}, so that there exists a rank variety 
$W^{r}_{E}(M)$ for $M\in \underline{\Mod}(kE)$ with the following properties:
\begin{itemize} 
\item[(i)]  $\widetilde W_{E}(M)\cong W_{E}^{r}(M)$, where $\widetilde W_E$ denotes the Benson-Carlson-Rickard extension \cite{BCR1} of the cohomological support 
for $E$ to $\underline{\Mod}(kE)$ and 
\item[(ii)] $W_{E}^{r}(M)= \varnothing$ if and only if $M$ is a projective $kE$-module. 
\end{itemize}
With the coproduct on $k[\Omega_{0}]$, we need to show that  
\begin{itemize} 
\item[(iii)]  $W_{E}^{r}(M\otimes N)=W_{E}^{r}(M)\cap W_{E}^{r}(N)$ for all $M, N \in \underline{\Mod}(k [\Omega_0])$ restricted to $E$-modules.
\end{itemize} 
Using the shifted subgroup description of the rank variety, the rank variety will consist of operators, 
$x_{\alpha}=\alpha_{1}t^{p^{r_{1}}-1}+\alpha_{2}t^{p^{r_{2}}-1}+\dots +\alpha_{s}t^{p_{r_{s}}-1}$ with $\alpha_{i}\in K$, $i=1,2,\dots, s$, 
where $K$ is a field extension of $k$ of large enough transcendence degree, such 
that $M_{\langle x_{\alpha} \rangle}$ is not free. 
We can now follow the argument in \cite[Lemma 3.9]{FP2}. The coproduct on $x_{\alpha}$ can be written as 
$$\Delta(x_{\alpha})\in x_{\alpha}\otimes 1 + 1\otimes x_{\alpha}+{\mathcal I}\otimes {\mathcal I},$$
where ${\mathcal I}$ is the augmentation ideal of $k[\Omega_0]$. By applying the argument in \cite[Proposition 2.2]{FP2}, it follows that $[M\otimes N]_{\langle \Delta(x_{\alpha}) \rangle}$ is free if and only if 
$[M\otimes N]_{\langle x_{\alpha}\otimes 1 + 1\otimes x_{\alpha} \rangle}$ is free. From this (iii) follows. 

Now we apply the fact that $\widetilde W_{k[\Omega_{0}]}(M)=\mathrm{res}^{*}(\widetilde W_{E}(M))$. By using the line of reasoning given in \cite[Theorems 10.6, 10.8]{BCR2}, 
one has 
\begin{itemize} 
\item[(iv)] $\widetilde W_{k[\Omega_{0}]}(M) = \varnothing$ if and only if $M$ is a projective $k[\Omega_{0}]$-module; 
\item[(v)]  $\widetilde W_{k[\Omega_{0}]}(M\otimes N)=\widetilde W_{k[\Omega_{0}]}(M)\cap \widetilde W_{k[\Omega_{0}]}(N)$. 
\end{itemize} 
This proves condition (c) of \thref{plavnik-witherspoon-classify} for $k[\Omega_{0}]$. Now \thref{plavnik-witherspoon-classify} yields the following:
\bpr{FGS-class} \cite[Theorem 10.3]{NP} For each finite group scheme $\Omega$,  
\[
\Spc \big( \underline{\modd}(k[\Omega]) \big) \mathrel{\mathop{\rightleftarrows}^{\rho}_\eta}\Proj \big( \opH^\bullet( k[\Omega_{0}], k) \big)^{\pi}
\]
are inverse homeomorphisms and the map
\[
\Phi_{W_C} : \ThickId \big(\underline{\modd}(k[\Omega]) \big) \to \mc{X}_{sp}\big(  \Proj ( \opH^\bullet( k[\Omega_{0}], k))^{\pi} \big)
\]
is an isomorphism of ordered monoids.
\epr
This result recovers a theorem of Negron and Pevtsova \cite[Theorem 10.3]{NP} who employed their hypersurface support theory. 
The example shows how it is obtained through a uniform approach based on Theorem D and \thref{plavnik-witherspoon-classify}. 

\appendix
\section{}
\label{App}

In this appendix we show that given any finite tensor category $\bT$ over a field $k$, 
there exists a compactly generated monoidal triangulated category $\widetilde{\bK}$ whose compact part is the stable category of $\bT$, that is, 
\[
\widetilde{\bK}^c \cong \underline{\bT}.
\]
We will construct $\widetilde{\bK}$ as the stable category of the indization of $\bT$. We recall this construction briefly. 

Ind-objects, and indizations of categories, were originally introduced in \cite{SGA}. We follow the construction as in \cite[Chapter 6]{KS1}. 
Consider the embedding of $\bT$ into $\bT^\vee$, 
the category of functors $\bT^{\op} \to \Set$, via the Yoneda embedding. 
Then the {\em{indization of $\bT$}} is the smallest full subcategory of $\bT^{\vee}$ containing the image of $\bT$, and closed under taking filtered colimits, and it is denoted $\Ind(\bT)$. 

We now claim the following, which is well-known to experts in the field. 
\bth{ind-objects}
Let $\bT$ be a finite tensor category over $k$. Then:
\begin{enumerate}
\item[(a)] $\Ind(\bT)$ is a Frobenius abelian monoidal category, and its stable category $\underline{\Ind(\bT)}$ is a compactly generated monoidal triangulated category.
\item[(b)] $\underline{\Ind(\bT)}^c \cong \underline{\bT}$. 
\end{enumerate}
\eth

\begin{proof}
Since $\bT$ is a finite tensor category over $k$, $\bT$ is equivalent (as an abelian category) to $\modd(A)$ for some finite dimensional $k$-algebra $A$ (see \cite[page 10]{EGNO}). Let $\Ind(\bT)$ be the ind-completion of $\bT$ as above. By \cite[Proposition 6.1.12]{KS1}, $\Ind(\bT) \cong \Mod(A)$, since $A$ is finite dimensional (and thus the finitely-presented $A$-modules are the finite dimensional $A$-modules). By \cite[Proposition 6.1.9, Proposition 6.1.12]{KS1}, a functor $F: \bT' \to \bT''$ extends to a functor $IF: \Ind(\bT') \to \Ind(\bT'')$, and $\Ind(\bT' \times \bT'') \cong \Ind(\bT') \times \Ind(\bT'')$, respectively. This implies that the tensor product $\otimes : \bT \times \bT \to \bT$ extends to a tensor product $\Ind(\bT) \times \Ind(\bT) \to \Ind(\bT)$. 

It is noted in \cite[Introduction]{KS1} that the category of ind-objects of a triangulated category does not appear to be triangulated; on the other hand, we claim that $\underline{\Ind(\bT)}$ {\em{is}} a triangulated category. Since $\bT$ is a finite tensor category, injectives and projectives coincide in $\bT$; in particular, this means that $A$ is injective over itself. By \cite[Theorem 15.9]{Lam1}, since $A$ is self-injective, if $M$ is any (not necessarily finite dimensional) $A$-module, $M$ is injective if and only if it is projective. Hence, injectives and projectives coincide for $\Ind(\bT)$ (i.e., it is Frobenius), and so the category $\underline{\Ind(\bT)}$ is a triangulated category. 

Now note that by \cite[Theorem 3]{Miyachi1}, if $R$ is a perfect ring, then $M$ is a compact object in $\underline{\Mod}(R)$ if and only if there is a finitely-generated $R$-module $M'$ with $M \cong M'$ in $\underline{\Mod}(R)$; since $A$ is finite dimensional, it is perfect, and so $\underline{ \Ind(\bT)}^c \cong \underline{\bT}$. 
\end{proof}


\section{}
\label{App2}

Let $\bK$ be a compactly generated M$\Delta$C with compact part $\bK^c$, equipped with the universal Balmer support $V$ sending objects of $\bK^c$ to closed sets in $\Spc \bK^c$. In this appendix, for $\bK$ such that $\Spc \bK^c$ is Noetherian, we construct an extension of $V$ to the non-compact objects of $\bK$.

We first recall a basic consequence of Noetherianity of $\Spc \bK^c$. 

\ble{desc-ideal}
Suppose $\Spc \bK^c$ is topologically Noetherian. Then finitely-generated thick ideals of $\bK^c$ satisfy the descending chain condition.
\ele

\begin{proof}
Suppose
\[
... \bI_i \subseteq \bI_{i-1} \subseteq ... \subseteq \bI_0
\]
is a descending chain of finitely-generated thick ideals of $\bK^c$. Then 
\[
... \Phi_V(\bI_i) \subseteq \Phi_V(\bI_{i-1}) \subseteq ... \subseteq \Phi_V(\bI_0)
\]
is a descending chain of closed subsets of $\Spc \bK^c$ (the fact that they are closed follows from finite generation, since if $\bI= \langle A_1,..., A_m \rangle$ then $\Phi_V(\bI)= \bigcup_{i=1...m} V(A_i)$). This chain stabilizes by the Noetherianity assumption. But note that every thick ideal is semiprime, e.g., equals the intersection of prime ideals over it (cf. \cite[Proposition 4.1.1]{NVY2}). This implies that if $\Phi_V(\bI_i) = \Phi_V(\bI_j)$, then $\bI_i = \bI_j$. 
\end{proof}

We next note a general lemma on localizing categories. Part (a) was originally proven in the symmetric case by Balmer--Favi \cite[Lemma 4.5]{BF1}.

\ble{int-loc}
We have the following.
\begin{enumerate}
\item[(a)] Suppose $\bI$ and $\bJ$ are thick ideals of $\bK^c$. Then $\Loc(\bI \cap \bJ)= \Loc(\bI) \cap \Loc(\bJ)$.
\item[(b)] Suppose $\Spc \bK^c$ is topologically Noetherian. Then for any set-indexed collection $\{\bI_i\}_{i \in I}$ of thick ideals of $\bK^c$, there is an equality $\Loc(\bigcap_{i\in I} \bI_i) = \bigcap_{i \in I} \Loc(\bI_i)$. 
\end{enumerate}
\ele

\begin{proof}
Part (a) follows exactly the same as in the symmetric case; we recall the Balmer--Favi proof here for reference. It is clear that $\Loc(\bI \cap \bJ) \subseteq \Loc(\bI) \cap \Loc(\bJ)$. By e.g., \cite[Theorem 4.3.3]{Neeman} or \cite[Theorem 7.2.1]{Krause}, 
\begin{align*}
\Loc(\bI) = &\{ A \in \bK : \text{for any morphism }C \xrightarrow{f} A \text{ with }C \in \bK^c,\text{ there exists } D \in \bI \\
&\text{ such that }f \text{ factors through }D\}
\end{align*}
for any thick subcategory $\bI$.
Let $A \in \Loc(\bI) \cap \Loc(\bJ)$, and suppose $C$ is compact with $f : C \to A$ a morphism. We must show that $f$ factors through some object in $\bI \cap \bJ$. By assumption, there exist objects $D_1 \in \bI$ and $D_2 \in \bJ$ with factorizations
\begin{center}
\begin{tikzcd}
C \arrow[rr, "f"] \arrow[rd] &                               & A             \\
                             & D_1 \arrow[ru] \arrow[r, "g"] & D_2 \arrow[u]
\end{tikzcd}
\end{center}
But now note that by standard properties of dual objects, $g: D_1 \to D_2$ factors 
\begin{center}
\begin{tikzcd}
D_1 \arrow[rrr, "g"] \arrow[rd, "\coev \otimes \id"] &                                                          &                                                             & D_2 \\
                                                     & D_1 \otimes D_1^* \otimes D_1 \arrow[r, "g \otimes \id"] & D_2 \otimes D_1^* \otimes D_1 \arrow[ru, "\id \otimes \ev"] &    
\end{tikzcd}
\end{center}
and hence the original morphism $f$ factors through $D_2 \otimes D_1^* \otimes D_1$, which is in $\bI \cap \bJ$. Therefore, $A \in \Loc(\bI \cap \bJ)$.

For (b), suppose $A$ is in $\bigcap_i \Loc(\bI_i)$, and again let $f: C \to A$ be a morphism with $C$ compact. Just as in part (a), we must show that $f$ factors through an object in $\bigcap_{i \in I} \bI_i$. Pick some $i \in I$ and set $\bI_i:=\bI_1$. We know by assumption that there is some $D_1 \in \bI_1$ such that $f$ factors through $D_1$. We now iterate this process: at the $k$th step, if for some $j$ we have $D_k \not \in \bI_j$, we can, using the same method as in part (a), construct $D_{k+1}$ such that $D_{k+1} \in \bI_j$ and $D_{k+1} \in \langle D_k \rangle$, and $f$ factors through $D_k$. This constructs a chain of thick ideals:
\[
... \langle D_{k+1} \rangle \subsetneq \langle D_{k} \rangle \subsetneq \langle D_{k-1} \rangle \subsetneq ... \subsetneq \langle D_1 \rangle.
\]
Now, since by assumption $\Spc \bK^c$ is Noetherian, it follows from \leref{desc-ideal} that this chain must terminate; in other words, at some step, the object $D:=D_i$ which has been constructed is in every ideal $\bI_j$ in our collection of ideals, since otherwise we could continue the chain of strict inclusions. By construction, then, $f$ factors through $D$, and $D \in \bigcap_i \bI_i$, which completes our proof. 

\end{proof}

We will assume for the remainder of this section that $\Spc \bK^c$ is Noetherian. We now define a version of extended Balmer supports for an arbitrary objects of $\bK$. Alternative methods for extending the Balmer support in the symmetric case include approaches via generalized tensor idempotents \cite{BF1, Stevenson1} based on the work of Rickard \cite{Rickard1}, via containment in homological primes \cite{Balmer3, BS1}, and via containment in ``big" primes \cite{BS1}. Another theory of supports for localizing ideals was developed in \cite{KL}, using the frame approach to localizing tensor ideals that was proposed in \cite{KP}.

\bde{ext-balmer-nc}
Let $A \in \bK$. Define the {\em{extended Balmer support of $A$}} as 
$$\widetilde V(A) = \{ \bP \in \Spc \bK^c : \forall C \in \bK^c \backslash \bP, A \otimes \bK^c \otimes C \not \subseteq \Loc(\bP)\}.$$
Thus $\widetilde V(-)$ forms a map from objects of $\bK$ to subsets of $\Spc \bK^c$.
\ede

\bth{props-nc}
If $C \in \bK^c$, then $\widetilde V(C) = V(C)$. Furthermore, $\widetilde V$ satisfies the following properties. In other words, $\widetilde V$ is an extended support which extends $V$. 
\begin{enumerate}
\item[(a)] $\widetilde V(0) = \varnothing$ and $\widetilde V(\unit)= \Spc \bK^c$.
\item[(b)] $\widetilde V(\bigoplus_{i \in I} A_i) = \bigcup_{i \in I} \widetilde V(A_i)$ for any set $I$ and collection of objects $A_i \in \bK$.
\item[(c)] If $A \to B \to D \to \Sigma A$ is a distinguished triangle, then $\widetilde V(A) \subseteq \widetilde V(B) \cup \widetilde V(D).$
\item[(d)] $\widetilde V(A) = \widetilde V(\Sigma A)$.
\item[(e)] $\bigcup_{D \in \bK^c} \widetilde V(A \otimes D \otimes C) = \widetilde V(A) \cap \widetilde V(C)$ for all $A \in \bK$ and $C \in \bK^c$. 
\end{enumerate}
\eth

\begin{proof}
First, let $C$ be compact. Then 
\begin{align*}
\widetilde V(C) &= \{\bP \in \Spc \bK^c : \forall D \in \bK^c \backslash \bP, \; C \otimes\bK^c \otimes D \not \subseteq \Loc(\bP) \}\\
&=\{\bP \in \Spc \bK^c : \forall D \in \bK^c \backslash \bP, \; C \otimes \bK^c \otimes D \not \in \bP\}\\
&=\{\bP \in \Spc \bK^c: C \not \in \bP\}\\
&= V(C).
\end{align*}
The second equality is by \cite[Lemma 2.2]{Neeman2} and the third is by the definition of a prime ideal. This shows that $\widetilde V$ extends the usual Balmer support defined on compact objects, and (a) follows from this since 0 and $\unit$ are both compact. 

Next, we show (b). We have $\bP \not \in \widetilde V(\bigoplus_{i \in I} A_i)$ if and only if there exists $C \in \bK^c \backslash \bP$ such that 
\[
\bigoplus_{i \in I} A_i \otimes \bK^c \otimes C \subseteq \Loc(\bP).
\]
Since $\Loc(\bP)$ is a localizing category, this implies that
\[
A_i \otimes \bK^c \otimes C \subseteq \Loc(\bP)
\]
for all $i$, and thus $\bP \not \in \widetilde V(A_i)$ for all $i$. This implies $\widetilde V(\bigoplus_{i \in I} A_i) \supseteq \bigcup_{i \in I} \widetilde V(A_i)$. 

For the other direction, suppose $\bP \not \in \widetilde V(A_i)$ for all $i$. Then for each $A_i$, there exists $C_i \in \bK^c \backslash \bP$ such that $A_i \otimes \bK^c \otimes C_i \subseteq \Loc(\bP)$. Consider the open set $\bigcup_{i \in I} \widetilde V(C_i)^c$ of $\Spc \bK^c$. By the assumption of Noetherianity, every subset of $\Spc \bK^c$ is quasi-compact, and so the cover of open sets $V(C_i)^c$ has a finite subcover, say $\widetilde V(C_1)^c,..., \widetilde V(C_n)^c$. 

By primeness of $\bP$, since each $C_i \not \in \bP$, we can pick compact objects $D_1,..., D_{n-1}$ such that $C:=C_1 \otimes D_1 \otimes C_2 \otimes D_2 \otimes... \otimes D_{n-1} \otimes C_n$ is not in $\bP$. We now claim that $C \in \langle C_i \rangle$ for all $i \in I$. To show this, we just need to show that $C$ is in every prime ideal $\bQ$ containing $C_i$, by \cite[Proposition 4.1.1]{NVY2}. Suppose $\bQ$ is a prime containing $\langle C_i \rangle$ for some $i \in I$. Then $\bQ \in \widetilde V(C_i)^c$, and hence in $\widetilde V(C_j)^c$ for some $j =1,..., n.$ In other words, $C_j \in \bQ$, and so it is now clear that $C \in \bQ$. Hence, $C \in \langle C_i \rangle$ for all $i \in I$.

But now note that for any $i \in I$, every object $D$ of $\langle C_i \rangle$ has the property that $A_i \otimes \bK^c  \otimes D \subseteq \Loc(\bP)$. This implies that for all $i$, we have $A_i \otimes \bK^c \otimes C \subseteq \Loc(\bP)$. Then
\[
(\bigoplus_{i \in I} A_i) \otimes \bK^c \otimes C \cong \bigoplus_{i \in I} (A_i \otimes \bK^c \otimes C) \subseteq \Loc(\bP),
\]
and, recalling that $C$ by construction was not in $\bP$, this shows that $\bP \not \in \widetilde V(\bigoplus_{i \in I} A_i)$. This completes the proof of (2).

For (c), suppose $A \to B \to D \to \Sigma A$ is a distinguished triangle in $\bK$, and that $\bP\not \in \widetilde V(B) \cup \widetilde V(D)$. Then there exist compact objects $C_1$ and $C_2$ which are not in $\bP$ such that $B \otimes \bK^c \otimes C_1$ and $D \otimes \bK^c \otimes C_2$ are both contained in $ \Loc(\bP)$. Since $C_1$ and $C_2$ are both not in $\bP$, there exists a compact object $E$ for which $C_1 \otimes E \otimes C_2 \not \in \bP$. But now it is clear that $A \otimes \bK^c \otimes C_1 \otimes E \otimes C_2 \subseteq \Loc(\bP)$, since $\Loc(\bP)$ is triangulated. Thus $\bP \not \in \widetilde V(A)$, and so $\widetilde V(A) \subseteq \widetilde V(B) \cup \widetilde V(D)$. 

(d) is straightforward from the fact that localizing categories are triangulated. 

For (e), suppose $\bP \not \in \bigcup_{D \in \bK^c} \widetilde V(A \otimes D \otimes C)$. Then for each $D \in \bK^c$, there exists $E_D \in \bK^c \backslash \bP$ such that $A \otimes D \otimes C \otimes \bK^c \otimes E_D \subseteq \Loc(\bP)$. Now, construct a sequence of ideals as follows. Begin with $\langle E_1 \rangle$ for some $E_{D_1}:=E_1$ corresponding to $D_1 \in \bK^c$ as above. Now suppose there exists an element $D_2$ in $\bK^c$ such that $E_2 \not \in  \langle E_1 \otimes F_1 \otimes E_2 \rangle$ for some $F_1$ such that $E_1 \otimes F_1 \otimes E_2 \not \in \bP$ (note that since $E_1$ and $E_2$ are both not in $\bP$ by assumption, it follows that $E_1 \otimes \bK^c \otimes E_2 \not \subseteq \bP$, and so there always exists at least one $F_1$ such that $E_1 \otimes F_1 \otimes E_2 \not \in \bP$). Then take $\langle E_1 \otimes F_1 \otimes E_2 \rangle$ as the next ideal, with $F_1$ as above. Continue in this manner. We obtain a sequence of strictly descending principal ideals generated by objects which are not in $\bP$:
\[
... \subsetneq \langle E_1 \otimes F_1 \otimes E_2 \otimes F_2 \otimes... \otimes E_i \rangle \subsetneq \langle E_1 \otimes F_1 \otimes... \otimes E_{i-1} \rangle \subsetneq ... \subsetneq \langle E_1 \rangle. 
\]
This sequence must terminate by Noetherianity and \leref{desc-ideal}. In other words, for some $n$, if we set $E:=E_1 \otimes F_1 \otimes... \otimes F_{n-1} \otimes E_n$, then for any $D \in \bK^c$ with corresponding object $E_D$, the object $E$ is in the ideal $\langle E \otimes F_D \otimes E_D \rangle$ for any $F_D$ such that $E \otimes F_D \otimes E_D \not \in \bP$. Now we claim that $A \otimes \bK^c \otimes C \otimes \bK^c \otimes E \subseteq \Loc(\bP)$. Given $D \in \bK^c$, we already know that $A \otimes D \otimes C \otimes \bK^c \otimes E_D \subseteq \Loc(\bP)$. But now note that any object $G$ in $\langle E \otimes F_D \otimes E_D \rangle$ also has the property that $A \otimes D \otimes C \otimes \bK^c \otimes G \subseteq \Loc(\bP)$, and in particular this is true for $E$, since $E \in \langle E \otimes F_D \otimes E_D\rangle$ by assumption. Therefore, $A \otimes D \otimes C \otimes \bK^c \otimes E \subseteq \Loc(\bP)$ for any $D \in \bK^c$, and thus $A \otimes \bK^c \otimes C \otimes \bK^c \otimes E \subseteq \Loc(\bP)$. 

The result now follows directly. If $\bP \in \widetilde V(C)$, then $C \otimes \bK^c \otimes E \not \subseteq \Loc(\bP)$, which means that there exists some $H \in \bK^c$ such that $C \otimes H \otimes E \not \in \bP$; and in that case, we can now see that $\bP \not \in \widetilde V(A)$, since this means that $A \otimes \bK^c \otimes (C \otimes H \otimes E) \subseteq \Loc(\bP)$. Therefore, $\bP \not \in V(A) \cap V(C)$, giving one direction of the claimed equality. 

For the other direction, suppose $\bP \in \bigcup_{D \in \bK^c} \widetilde V(A \otimes D \otimes C)$. Then there is some $D$ for which $\bP \in \widetilde V(A \otimes D \otimes C)$. Suppose $\bP \not \in \widetilde V(C)$, in other words, $C \in \bP$. Then $A \otimes D \otimes C$ would be in $\Loc(\bP)$, which would immediately imply $\bP \not \in \widetilde V(A \otimes D \otimes C)$, a contradiction. In other words, $\bP \in \widetilde V(C)$. We must just show now that $\bP \in \widetilde V(A)$. Suppose $E \in \bK^c \backslash \bP$. If $A \otimes \bK^c \otimes E \subseteq \Loc(\bP)$, then this would imply that $A \otimes D \otimes C \otimes \bK^c \otimes E \subseteq \Loc(\bP)$, which would contradict the fact that $\bP \in \widetilde V(A \otimes D \otimes C)$. Hence $A \otimes \bK^c \otimes E \not \subseteq \Loc(\bP)$ for all $E \in \bK^c \backslash \bP$, which implies $\bP \in \widetilde V(A)$. This completes the proof. 
\end{proof}

Next we give a characterization of elements with empty support for categories with finite Krull dimension. For ease of notation, for any integer $m \geq 1$, we will denote $A(m)$ for the collection of objects
\[
A(m):=\underbrace{A \otimes \bK^c \otimes A \otimes ... \otimes \bK^c \otimes A}_\text{$m$ copies of $A$}.
\]

\bth{faith-nilp-nc}
Suppose $A \in \bK$ satisfies $\widetilde V(A)= \varnothing$. If $\bK^c$ has finite Krull dimension $n$ (that is, any chain of strict containments of prime ideals is at most length $n$), then $A(n)=\{0\}$. 
\eth

\begin{proof}
Suppose $\widetilde V(A)= \varnothing$. Then for each prime $\bP \in \Spc \bK^c$, there exists a compact object $C$ which is not in $\bP$, such that $A \otimes \bK^c \otimes C \subseteq \Loc(\bP)$. Since $\bigcap_{\bP \in \Spc \bK^c} \bP = \{0\}$, it is enough to show that $A(n) \subseteq \Loc(\bP)$ for all primes $\bP$, by \leref{int-loc}. 

We will show by induction on the maximum $m$ such that there exists a chain of prime ideals $\bP= \bP_1 \subsetneq \bP_2 \subsetneq ... \subsetneq \bP_m$ that $A(m) \subseteq \Loc(\bP)$. If $\bP$ is a maximal ideal, then $\langle \bP, C \rangle = \bK^c$ for any $C \in \bK^c \backslash \bP$, in particular for some $C$ for which $A \otimes \bK^c \otimes C \subseteq \Loc(\bP)$, which exists by assumption. Hence one can form the unit object $\unit$ by successively taking shifts, cones, direct sums, and tensor products with arbitrary compact objects, applied to objects in $\{\bP, C\}$. Note that each operation preserves the property that applying $A \otimes \bK^c \otimes -$ gives a collection of objects in $\Loc(\bP)$; since objects of $\bP$, as well as $C$ by assumption have this property, so does $\unit$. In other words, $A \otimes \bK^c \otimes \unit \subseteq \Loc(\bP)$, and in particular $A(1) =\{A\} \subseteq \Loc(\bP)$. This completes the base case of the induction.

Now, assume that $\bP$ is any prime in $\Spc \bK^c$, with $m$ the maximum number such that there is a chain of strict containments of prime ideals over $\bP$ of length $m$. By assumption, as before, there exists a compact $C$ not in $\bP$ with $A \otimes \bK^c \otimes C \subseteq \Loc(\bP)$. Now note that it is enough to show that $A(m-1) \subseteq \Loc \langle \bP, C \rangle$; if this is true, one can form any object in this set by successively taking shifts, cones, arbitrary set-indexed direct sums, and tensor products, applied to objects in $\{\bP, C\}$. As in the base case of the induction, each of these steps preserves the property that applying $A \otimes \bK^c \otimes -$ produces a collection of objects contained in $\Loc(\bP)$. Therefore, this implies that $A(m) \subseteq \Loc(\bP)$. 

To show that $A(m-1)$ is contained in $\Loc\langle \bP, C\rangle$, it is enough to show that it is contained in $\Loc(\bQ)$ for any prime ideal $\bQ$ containing $\langle \bP, C \rangle$, since every thick ideal is the intersection of prime ideals above it (again, by \cite[Proposition 4.1.1]{NVY2}) and by \leref{int-loc}. But such a prime ideal $\bQ$ properly contains $\bP$, and so by the inductive hypothesis $A(m-1)$ is contained in $\Loc(\bQ)$ for all such $\bQ$. This completes the proof. 
\end{proof}

\noindent 
{\it On behalf of all authors, the corresponding author states that there is no conflict of interest.}

\end{document}